\documentclass[12pt,oneside]{amsart}
\usepackage{latexsym}
\usepackage{amsmath,amssymb,amsthm}
\usepackage{amscd}
\usepackage[dvips]{graphics}

\newcommand{\re}{\mathbb{R}}
\newcommand{\co}{\mathbb{C}}

\newcommand{\Sl}[1]{\mathbf{SL}(#1,\mathbb{R})}

\newcommand{\na}{\nabla}
\newcommand{\sfrac}[2]{{\textstyle \frac{#1}{#2}}}
\newcommand{\D}{\displaystyle}

\newtheorem{prop}{Proposition} %[subsection]
\newtheorem{cor}[prop]{Corollary}
\newtheorem{thm}{Theorem}

\newtheorem{lem}[prop]{Lemma}

\theoremstyle{remark}
\newtheorem*{rem}{Remark}

\begin{document}
\author{John Loftin}
\title{Affine Hermitian-Einstein Metrics}
\maketitle

\section{Introduction}

A holomorphic vector bundle $E\to N$ over a compact K\"ahler
manifold $(N,\omega)$ is called \emph{stable} if every coherent
holomorphic subsheaf $F$ of $E$ satisfies $$ 0 < \mbox{rank } F  <
\mbox {rank } E \qquad \Longrightarrow \qquad
\mu_\omega(F)<\mu_\omega (E),$$ where $\mu_\omega$ is the
$\omega$-\emph{slope} of the sheaf given by $$\mu_\omega(E) =
\frac{\deg_\omega(E)}{\mbox{rank }E} = \frac{\int_N c_1(E,h) \wedge
\omega^{n-1}}{\mbox{rank }E}.$$ Here $c_1(E,h)$ is the first Chern
form of $E$ with respect to a Hermitian metric $h$.  The famous
theorem of Donaldson \cite{donaldson85, donaldson87} (for algebraic
manifolds only) and Uhlenbeck-Yau \cite{uhlenbeck-yau86,
uhlenbeck-yau89} says that an irreducible vector bundle $E\to N$ is
$\omega$-stable if and only if it admits a Hermitian-Einstein metric
(i.e.\ a metric whose curvature, when the 2-form part is contracted
with the metric on $N$, is a constant times the identity
endomorphism on $E$).

An important generalization of this theorem is provided by
Li-Yau \cite{li-yau87} for complex manifolds (and subsequently
due to Buchdahl by a different method for surfaces
\cite{buchdahl88}). The major insight for this extension is the
fact that the degree is well-defined as long as the Hermitian
form $\omega$ on $N$ satisfies only $\partial\bar\partial
\omega^{n-1} =0$. This is because
$$\deg_\omega (E) = \int_N c_1(E,h)\wedge \omega^{n-1}$$ and the
difference of any two first Chern forms $c_1(E,h)-c_1(E,h')$ is
$\partial\bar\partial$ of a function on $N$.  But then Gauduchon has
shown that such an $\omega$ exists in the conformal class of every
Hermitian metric on $N$ \cite{gauduchon75,gauduchon77}.  (Such a
metric on $N$ is thus called a Gauduchon metric.)  The book of
L\"ubke-Teleman \cite{lubke-teleman95} is quite useful, in that it
contains most of the relevant theory in one place.

An affine manifold is a real manifold $M$ which admits a flat,
torsion-free connection $D$ on its tangent bundle. It is well
known (see e.g.\ \cite{shima07}) that $M$ is an affine manifold
if and only if $M$ admits an affine atlas whose transition
functions are locally constant elements of the affine group
$$\mbox{Aff}(n)=\{\Phi\!: \re^n \to \re^n, \quad \Phi\!:x\mapsto
Ax+b\}.$$ (In this case, geodesics of $D$ are straight line segments
in the coordinate patches of $M$.) The tangent bundle $TM$ of an
affine manifold admits a natural complex structure, and it is often
fruitful to think of $M$ as a real slice of a complex manifold. In
particular, local coordinates $x=(x^1,\dots,x^n)$ on $M$ induce the
local frame $y=(y^1,\dots,y^n)$ on $TM$ so that every tangent vector
$y$ can be written as $y=y^i\frac\partial{\partial x^i}$.  Then
$z^i=x^i+\sqrt{-1}y^i$ form holomorphic coordinates on $TM$.  We
will usually denote the complex manifold $TM$ as $M^\co$.

Cheng-Yau \cite{cheng-yau82} proved the existence of affine
K\"ahler-Einstein metrics on appropriate affine flat manifolds.  The
setting in this case is that of affine K\"ahler, or Hessian, metrics
(see also Delano\"e \cite{delanoe89} for related results).  A
Riemannian metric $g$ on $M$ is \emph{affine K\"ahler} if each point
has a neighborhood on which there are affine coordinates $\{x^i\}$
and a real potential function $\phi$ satisfying $$ g_{ij}\,dx^i dx^j
= \frac{\partial^2\phi}{\partial x^i
\partial x^j} \,dx^i dx^j.$$ Every Riemannian metric $g$ on $M$
extends to a Hermitian metric $g_{ij} \,dz^i \overline {dz^j}$ on
$TM$.  The induced metric on $M^\co$ is K\"ahler if and only if the
original metric is affine K\"ahler.

An important class of affine manifolds is the class of \emph{special
affine manifolds}, those which admit a $D$-covariant constant volume
form $\nu$. If such an affine manifold admits an affine K\"ahler
metric, then Cheng-Yau showed that the metric can be deformed to a
flat metric by adding the Hessian of a smooth function
\cite{cheng-yau82}.  There is also the famous conjecture of Markus:
A compact affine manifold admits a covariant-constant volume form if
and only if $D$ is complete. In the present work, we will use a
covariant-constant volume form to convert $2n$-forms on the complex
manifold $TM=M^\co$ to $n$-forms on $M$ which can be integrated. The
fact that $D\nu=0$  will ensure that $\nu$ does not provide
additional curvature terms when integrating by parts on $M$.

The correct analog of a holomorphic vector bundle over a complex
manifold is a flat vector bundle over an affine manifold.  In
particular, the transition functions of a real vector bundle over an
affine flat $M$ may be extended to transition functions on $TM$ by
making them constant along the fibers of $M^\co\to M$.  In the local
coordinates as above, we require the transition functions to be
constant in the $y$ variables.  Such a transition function $f$ is
holomorphic over $TM$ exactly when $$0=\bar\partial f = \frac
{\partial f} {\partial \overline{z^i}}\,\overline{dz^i} = \left(
\frac12 \frac{\partial f}{\partial x^i} + \frac{\sqrt{-1}}2
\frac{\partial f}{\partial y^i} \right)\overline{dz^i} = \frac12
\frac{\partial f}{\partial x^i} \,\overline{dz^i},$$ in other words,
when the transition function is constant in $x$. In this way, from
any locally constant vector bundle $E\to M$, we can produce a
locally constant holomorphic vector bundle of the same rank
$E^\co\to M^\co$.

The existence of Hermitian-Einstein metrics on holomorphic
vector bundles over Gauduchon surfaces has been used by
Li-Yau-Zheng \cite{li-yau-zheng90,li-yau-zheng94} (and also
Teleman \cite{teleman94}, based on ideas in
\cite{li-yau-zheng90}) to provide a new proof of Bogomolov's
theorem on compact complex surfaces in Kodaira's class VII$_0$.

The theory we present below is explicitly modeled on Uhlenbeck-Yau
and Li-Yau's arguments.  We have found it useful to follow the
treatment of L\"ubke-Teleman \cite{lubke-teleman95} fairly closely,
since most of the relevant theory for Hermitian-Einstein metrics on
Gauduchon manifolds is contained in \cite{lubke-teleman95}.  Our
main theorem is

\begin{thm} \label{main-thm}
Let $M$ be a compact special affine manifold without boundary
equipped with an affine Gauduchon metric $g$.  Let $E\to M$ be a
flat complex vector bundle. If $E$ is $g$-stable, then there is an
affine Hermitian-Einstein metric on $E$.
\end{thm}

A similar result holds for flat real vector bundles over $M$ (see
Corollary \ref{re-cor} below).

We should remark that the affine K\"ahler-Einstein metrics
produced by Cheng-Yau in \cite{cheng-yau82} are examples of
affine Hermitian-Einstein metrics as well: The affine
K\"ahler-Einstein metric $g$ on the affine manifold $M$ can be
thought of as a metric on the flat vector bundle $TM$, and as
such a bundle metric, $g$ is affine Hermitian-Einstein with
respect to $g$ itself as an affine K\"ahler metric on $M$.
Cheng-Yau's method of proof is to solve real Monge-Amp\`ere
equations on affine manifolds (and they also provide one of the
first solutions to the real Monge-Amp\`ere equation on convex
domains in \cite{cheng-yau82}).

It is worth pointing out, in broad strokes, how to relate the proof
we present below to the complex case:  The complex case relies on
most of the standard tools of elliptic theory on compact manifolds:
the maximum principle, integration by parts, $L^p$ estimates,
Sobolev embedding, spectral theory of elliptic operators, and some
intricate local calculations.  The main innovation we provide to the
affine case is Proposition \ref{int-by-parts} below, which secures
our ability to integrate by parts on a special affine manifold.
Moreover, by extending a complex flat vector bundle $E\to M$ to a
flat holomorphic vector bundle $E^\co\to M^\co$ as above, we can
ensure that the local calculations on $M$ are \emph{exactly the
same} as those on $M^\co$, and thus we do not have to change these
calculations at all to use them in our proof.  The maximum principle
and spectral theory work the same way in our setting as well. The
$L^p$ and Sobolev theories in the complex case do not strongly use
the ambient real dimension $2n$ of the complex manifold:  and in
fact, reducing the dimension to $n$ helps matters.

There are a few other small differences in our approach on affine
manifolds as compared to the case of complex manifolds:  First of
all, we are able to avoid the intricate proof of Uhlenbeck-Yau
\cite{uhlenbeck-yau86,uhlenbeck-yau89} that a weakly holomorphic
subbundle of a holomorphic vector bundle on a complex manifold is a
reflexive analytic subsheaf (see also Popovici \cite{popovici05}).
The corresponding fact we must prove is that a weakly flat subbundle
of a flat vector bundle on an affine manifold is in fact a flat
subbundle.  We are able to give a quite simple regularity proof in
the affine case below in Proposition \ref{produce-pi}, and the flat
subbundle we produce is smooth.

Another small difference  between the present case and the complex
case concerns simple bundles. The important estimate Proposition
\ref{orthog-kernel} below works only for simple bundles $E$ (bundles
whose only endomorphisms are multiples of the identity). This does
not affect the main theorem in the complex case, for Kobayashi
\cite{kobayashi72} has shown that any stable holomorphic vector
bundle over a compact Gauduchon manifold must be simple.  For a flat
\emph{real} vector bundle $E$ over an affine manifold, there are two
possible notions of simple, depending on whether we require every
real locally constant section of $\mbox{End}(E)$ ($\re$-simple), or
every complex locally constant section of
$\mbox{End}(E)\otimes_\re\co$ ($\co$-simple), to be a multiple of
the identity. Since Kobayashi's proof relies on taking an
eigenvalue, we must do a little more work in Section
\ref{simple-section} below to address the case of $\re$-simple
bundles.

In Sections \ref{aff-dol-sec} and \ref{flat-vec-sec} below, we
develop some of the basic theory of $(p,q)$ forms with values in a
flat vector bundle $E$ over $M$, affine Hermitian connections, and
the second fundamental form. The basic principle behind these
definitions is to mimic the same formulas of the holomorphic vector
bundle $E^\co\to M^\co$. One interesting side note in this story is
Lemma \ref{dual-Hermitian}, which notes for a metric on a real flat
vector bundle $(E,\na)$ over $M$, the dual connection $\na^*$ on $E$
is equivalent to the Hermitian connection on $E^\co\to M^\co$.

 Section \ref{int-by-parts-sec} contains our
main technical tool, which allows us to integrate $(p,q)$ forms
by parts on a special affine manifold.  Then in Section
\ref{ahem-sec}, we prove the easy parts of the theory of affine
Hermitian-Einstein metrics: vanishing, uniqueness, and
stability theorems for affine Hermitian-Einstein metrics, most
of which are due to Kobayashi in the complex case. The proofs
we present are easier than in the complex case, since we need
only consider subbundles, and not singular subsheaves, in our
definition of stability. In Section \ref{aff-g-section}, we
produce affine Gauduchon metrics on special affine manifolds.

Then in Sections \ref{cont-meth-sec} to \ref{dest-sub-sec}, we
prove Theorem \ref{main-thm}, following the continuity method
of Uhlenbeck-Yau, as modified by Li-Yau for Gauduchon manifolds
and as presented in L\"ubke-Teleman \cite{lubke-teleman95}.
Since our local calculations are designed to be exactly the
same as the complex case, we omit some of these calculations.
On the other hand, we do emphasize those parts of the proof
which involve integration, as this highlights the main
difference between our theory on affine manifolds and the
complex case. The regularity result in Section
\ref{dest-sub-sec} is much easier than that of Uhlenbeck-Yau
\cite{uhlenbeck-yau86, uhlenbeck-yau89}. Finally, in Section
\ref{simple-section}, we address the issue of $\re$- and
$\co$-simple bundles, to prove a version of the main theorem,
Corollary \ref{re-cor} for $\re$-stable flat real vector
bundles.

We should also mention Corlette's results on flat principle bundles
on Riemannian manifolds:
\begin{thm} \cite{corlette88}
Let $G$ be a semisimple Lie group,  $(M,g)$ a compact Riemannian
manifold, and $P$ a flat principle $G$-bundle over $M$. A
\emph{metric} on $P$ is defined to be a reduction of the structure
group to $K$ a maximal compact subgroup of $G$, and a \emph{harmonic
metric} is a metric on $P$ so that the induced
$\pi_1(M)$-equivariant map from the universal cover $\tilde M$ to
the Riemannian symmetric space $G/K$ is harmonic. Then $P$ admits a
harmonic metric if and only if $P$ is reductive in the sense that
the Zariski closure of the holonomy at every point in $M$ is a
reductive subgroup of $G$.
\end{thm}

If $G$ is the special linear group, then we may consider the flat
vector bundle $(E,\na)$ associated to $P$. Then the reductiveness of
$P$ is equivalent to the condition on $E$ that any $\na$-invariant
subbundle has a $\na$-invariant complement.  For $M$ a compact
special affine manifold equipped with an affine Gauduchon metric
$g$, our Theorem \ref{main-thm} produces an affine
Hermitian-Einstein bundle metric on a flat vector bundle $E$ when it
is slope-stable.  If we assume $E$ is irreducible as a flat bundle,
then our slope-stability condition is \emph{a priori} weaker than
Corlette's: we require every proper flat subbundle of $E$ to have
smaller slope, while Corlette requires that there be no proper flat
subbundles of $E$.  It should be interesting to compare the harmonic
and affine Hermitian-Einstein metrics on $E$ when they both exist.

It is well known that an affine structure on a manifold $M$ is
equivalent to the existence of an \emph{affine-flat} (flat and
torsion-free) connection $D$ on the tangent bundle $TM$, which
induces a flat connection on a principle bundle over $M$ with group
$G= {\rm Aff}(n,\re)$ the affine group. The affine group is not
semisimple (or even reductive), and so Corlette's result does not
apply directly to study this case.  On a special affine manifold,
however, $D$ induces a flat metric on an $\sl n$-principal bundle,
and Corlette's result applies on $TM$ as a flat $\sl n$-bundle.
Thus, Corlette's result cannot see that $D$ is torsion-free. On the
other hand, the affine Hermitian-Einstein metric we produce does
essentially use the fact that $D$ is torsion-free: this ensures the
induced almost-complex structure on $M^\co$ is integrable.  So the
affine Hermitian-Einstein metrics should be able to exploit the
affine structure on $M$.

I would like to thank D.H.\ Phong, Jacob Sturm, Bill Goldman
and S.T.\ Yau for inspiring discussions.  I am also grateful to
the NSF for support under grant DMS0405873.

\section{Affine Dolbeault complex} \label{aff-dol-sec}

On an affine manifold $M$, there are natural $(p,q)$ forms (see
Cheng-Yau  \cite{cheng-yau82} or Shima \cite{shima07}), which
are the restrictions of $(p,q)$ forms from $M^\co$. We define
the space of these forms as
$$\mathcal A^{p,q}(M) = \Lambda^p(M) \otimes \Lambda^q(M)$$
for $\Lambda^p(M)$ the usual exterior $p$-forms on $M$. If $x^i$ are
local affine coordinates on $M$, then we will denote the induced
frame on $\mathcal A^{p,q}$ by $$\{ dz^{i_1} \wedge \cdots \wedge
dz^{i_p} \otimes d\bar z^{j_1} \wedge \cdots \wedge d\bar z^{j_q}
\},$$ where we think of $z^i=x^i+\sqrt{-1} \,y^i$ as coordinates on
$M^\co$ as above.

The flat connection $D$ induces flat connections on the bundles
$\Lambda^q(M)$ of $q$-forms of $M$. Therefore, the exterior
derivative $d$ extends to operators
$$\begin{array}{c} d^D\otimes I\!: \Lambda^p(M) \otimes \Lambda^q(M)\to
\Lambda^{p+1}(M) \otimes \Lambda^q(M), \\[2mm]
 I\otimes d^D\!:
\Lambda^p(M) \otimes \Lambda^q(M)\to \Lambda^p(M) \otimes
\Lambda^{q+1}(M).
\end{array}$$
 for $I$ the identity operator and $d^D$ the exterior derivative
for bundle-valued forms induced from $D$.  These operators are
equivalent to the operators $\partial$ and $\bar\partial$ restricted
from $M^\co$.  We find it useful to use the exact restrictions of
$\partial$ and $\bar \partial$ (so that, insofar as possible, all
the local calculations we do are the same as in the case of complex
manifolds).  The proper correspondences are, for $\partial$ and
$\bar \partial$ acting on $(p,q)$ forms,
$$ \partial = \sfrac12 (d^D \otimes I), \qquad \bar \partial = (-1)^p
\sfrac12 (I \otimes d^D).$$

A Riemannian metric $g$ on $M$ gives rise to a natural $(1,1)$
form given in local coordinates by $\omega_g =
g_{ij}dz^i\otimes d\bar z^j$. This is of course the restriction
of the Hermitian form induced by the extension of $g$ to
$M^\co$.

There is also a natural wedge product on $\mathcal A^{p,q}$, which
we take to be the restriction of the wedge product on $M^\co$: If
$\phi_i \otimes \psi_i\in \mathcal A^{p_i,q_i}$ for $i=1,2$, then we
define
$$ (\phi_1\otimes \psi_1) \wedge (\phi_2\otimes \psi_2) =
(-1)^{q_1p_2} (\phi_1\wedge\phi_2)\otimes(\psi_1\wedge \psi_2) \in
\mathcal A^{p_1+p_2, q_1+q_2}.$$

Consider the space of $(p,q)$ forms $\mathcal A^{p,q}(E)$ taking
values in a complex (or real) vector bundle $E\to M$.  If $\na$ is a
flat connection on $E$, and $h$ is a Hermitian metric on $E$
(positive-definite if $E$ is a real bundle), then we consider the
Hermitian connection, or Chern connection, on $E^\co\to M^\co$.
Recall the Hermitian connection is the unique connection on a
Hermitian holomorphic vector bundle over a complex manifold which
both preserves the Hermitian metric and whose $(0,1)$ part is equal
to the natural $\bar\partial$ operator on sections of $E$. Any
locally constant frame $s_1,\dots,s_r$ over $E\to M$ extends to a
holomorphic frame over $E^\co\to M^\co$, where we have the usual
formula (see e.g.\ \cite{kobayashi87}) for the Hermitian connection:
If $h_{\alpha\bar\beta}= h(s_\alpha,s_\beta)$, then the connection
form is a  $\mbox{End} E$-valued $(1,0)$ form
$$\theta^\alpha_\beta = h^{\alpha\bar \gamma} \,\partial
h_{\beta\bar \gamma}.$$

In passing from $(p,q)$ forms on $M^\co$ to $(p,q)$ forms on $M$, we
use the following convention:
\begin{equation}
 \label{form-convention}
 dz^{i_1}\wedge \cdots \wedge dz^{i_p} \wedge d\bar z^{j_1} \wedge
 \cdots \wedge d\bar z^{j_q} \mapsto dz^{i_1}\wedge \cdots \wedge dz^{i_p} \otimes d\bar z^{j_1} \wedge
 \cdots \wedge d\bar z^{j_q}
\end{equation}
As we will see in the next section, this convention will make all
the important curvature quantities on $E\to M$ to be real in the
case $E$ is a real vector bundle equipped with a real
positive-definite metric.

There is also a natural map from $(p,q)$ forms on $M$ to $(q,p)$
forms on $M$, which is the restriction of complex conjugation on
$M^\co$: If $\alpha\in \Lambda^p(M)$, $\beta\in \Lambda^q(M)$ are
complex valued forms, then we define
\begin{equation} \label{conj-form}
\overline{\alpha\otimes \beta} = (-1)^{pq} \bar \beta \otimes \bar
\alpha.
\end{equation}

At least when $E$ is a real bundle and $h$ is a real
positive-definite metric, the Hermitian connection described above,
when restricted to $M$, has an interpretation in terms of the dual
connection of $\na$ with respect to $h$. Recall that the dual
connection $\na^*$ is defined on $E\to M$ by
$$ d[ h(s_1,s_2)] = h(\na s_1, s_2) + h(s_1, \na^* s_2)$$ (see e.g.\
\cite{amari-nagaoka}). Then we may define operators
$\partial^{\na,h}$ and $\bar\partial^\na$ on $\mathcal A^{p,q}(E)$
as follows: For $\phi\in \mathcal A^{p,q}$ and $\sigma\in
\Gamma(E)$,
\begin{eqnarray*}
 \partial^{\na,h} \sigma &=& \na^* \sigma \otimes \sfrac12, \\
 \bar\partial^\na \sigma &=& \sfrac12 \otimes \na\sigma, \\
 \partial^{\na,h}(\sigma\cdot \phi) &=& (\partial^{\na,h}\sigma) \wedge \phi +\sigma
 \cdot \partial \phi, \\
\bar\partial^\na(\sigma\cdot \phi) &=& (\bar\partial^\na\sigma)
\wedge \phi +\sigma \cdot \bar\partial \phi.
\end{eqnarray*}
On $M$, we consider the pair $(\partial^{\na,h},\bar\partial^\na)$
to form an \emph{extended Hermitian connection} on $E$, and the
extended connection is equivalent to the Hermitian connection on
$E^\co\to M^\co$: The Hermitian connection on $E^\co\to M^\co$ is
given by $d^{\na,h} =
\partial^{\na,h} + \bar\partial^{\na} \!: \Lambda^0(E^\co) \to
\Lambda^1(E^\co)$.

 Also note that the difference $\na^*-\na$ is a
section of $\Lambda^1(\mbox{End} E)$.  This is a similar
construction to the first Koszul form on a Hessian manifold (see
e.g.\ Shima \cite{shima07}).

We have the following lemma, whose proof is a simple computation:
\begin{lem} \label{dual-Hermitian}
If $(E,\na)$ is a flat real vector bundle over an affine manifold
$M$, and $E$ is equipped with a positive-definite metric $h$, then
the extended Hermitian connection on $E$ (when considered as a
complex vector bundle with Hermitian metric induced from $h$) is
given by $$(\partial^{\na,h}, \bar\partial^\na) = (\na^* \otimes
\sfrac12 , \sfrac12 \otimes \na)$$ for $\na^*$ the dual connection
of $\na$ on $E$ with respect to the metric $h$.
\end{lem}

The curvature form $\Omega\in\mathcal A^{1,1}(\mbox{End} E)$ is
given by $$\Omega^\alpha_\beta = \bar
\partial \theta^\alpha_\beta =  -h^{\alpha\bar\eta}\partial\bar\partial h_{\beta\bar\eta} +
 h^{\alpha\bar\zeta} h^{\epsilon\bar\eta} \partial h_{\beta\bar\eta}
 \wedge \bar\partial h_{\epsilon\bar\zeta}.$$
If we write $\Omega^\alpha_\beta = R^\alpha_{\beta i\bar\jmath}
\,  dz^i \wedge d\bar z^j,$ then
$$R^\alpha_{\beta i\bar\jmath} = -h^{\alpha\bar\eta}
\frac{\partial^2 h_{\beta\bar\eta}} {\partial z^i \partial \bar z^j}
 + h^{\alpha\bar \zeta} h^{\epsilon\bar \eta} \frac{\partial h_{\beta
 \bar\eta}} { \partial z^i} \frac{\partial h_{\epsilon \bar \zeta}}
 {\partial \bar z^j}.$$
These same formulas represent the restriction of the curvature
form of $E^\co\to M^\co$ to $M$.  On $M$, we call this the
\emph{extended curvature form} (and we still use the symbols
$dz^i,d\bar z^j$ to represent elements of $\mathcal A^{p,q}$ on
$M$).

We use a Riemannian metric $g$ on $M$ to contract the $(1,1)$
part of an extended curvature form to form a section of
$\mbox{End} E = E^*\otimes E$ which we call the \emph{extended
mean curvature}. A metric on $E$ is said to be \emph{affine
Hermitian-Einstein} with respect to $g$ if its extended mean
curvature $K^\alpha_\beta$ is a constant $\gamma$ times the
identity endomorphism of $E$.  In index notation, we have
$$K^\alpha_\beta = g^{i\bar\jmath} R^\alpha_{\beta i\bar\jmath}
= \gamma \, I^\alpha_\beta.$$ (Here we extend the Riemannian
metric $g$ to a Hermitian metric $g_{i\bar\jmath}$ on $M^\co$,
and $I$ is the identity endomorphism on $E$.)

Given a Hermitian locally constant bundle $(E,h)$ on $M$, the trace
$R^\alpha_{\alpha i\bar\jmath}$ is called the \emph{extended first
Chern form}, or extended Ricci curvature. This first Chern form is
give by
$$c_1(E,h) = -\partial \bar \partial \log\det h_{\alpha\bar
\beta},$$ and it may naturally be thought of as the extended
curvature of the locally constant line bundle $\det E$ with metric
$\det h$.

The extended first Chern form and the extended mean curvature are
related by
 \begin{equation} \label{trace-K}
 (\mbox{tr}\, K)\, \omega_g^n =n\, c_1(E,h)\wedge \omega_g^{n-1}.
\end{equation}

\section{Flat vector bundles} \label{flat-vec-sec}
In this section, we collect some facts about flat vector bundles,
and representations of the fundamental group, and vector-bundle
second fundamental forms. The field $\mathbb K$ will represent
either $\re$ or $\co$.

A section $s$ of a flat vector bundle $(E,\na)$ over a manifold
$M$ is called \emph{locally constant} if $\na s=0$.  Every flat
vector bundle has local frames of locally constant sections,
given by parallel transport from a basis of a fiber $E_x$ for
$x\in M$.  For this reason, flat vector bundles are sometimes
referred to as locally constant vector bundles.

A flat $\mathbb K$-vector bundle of rank $r$ naturally corresponds
to a representation $\rho$ of fundamental group into $\mathbf{GL}
(r,\mathbb K)$. For $\tilde M$ the universal cover of $M$, the
fundamental group $\pi_1(M)$ acts on total space $\tilde M \times
\mathbb K^r$  equivariantly with respect to the action
$$\gamma \!: (x,y) \to (\gamma(x), \rho(\gamma)(y)).$$

In this picture, a flat subbundle of rank $r'$ is given by an
inclusion $\tilde M\times \mathbb K^{r'} \subset \tilde M \times
\mathbb K^r$ as trivial bundles, where $\pi_1$ acts on $\tilde
M\times \mathbb K^{r'}$. In other words, we require for every
$\gamma\in\pi_1$ and $y\in \re^{r'}$, $\rho(\gamma)(y) \in
\re^{r'}$.

Let $(E,\na)$ be a flat complex vector bundle over an affine
manifold $M$, and $h$ is a Hermitian metric on $E$.  The
geometry of flat subbundles of $E$ follows as in the case of
holomorphic bundles on complex manifolds.  Let $F$ be a flat
subbundle of $E$ (i.e. $F$ is a smooth subbundle of $E$ whose
sections $s$ satisfy $\na_X s$ is again a section of $F$ for
every vector field $X$ on $M$). Then for any section $s$ of
$F$, we may split $\partial^{\na,h}s$ into a part in $F$ and a
part $h$-orthogonal to $F$:
$$\partial^{\na,h} s =
\partial^{\na_F, h_F} s + A(s).$$ As the notation suggests, the
first term on the right $\partial^{\na_F, h_F} s$ is the $(1,0)$
part of the affine Hermitian connection induced on $F$ by $\na$ and
$h$. The second term $A$ is a ${\rm Hom}(F,F^\perp)$-valued $(1,0)$
form called the \emph{second fundamental form} of the subbundle $F$
of $E$. Note we only need consider $\partial^{\na,h}s$ since the
second fundamental form is of $(1,0)$ type in the complex case. We
have the following  proposition (see e.g. \cite[Proposition
I.6.14]{kobayashi87})
\begin{prop} \label{second-fund-prop}
Given $(E,\na)$, $h$, $F$ and $A$ as above, $A$ vanishes identically
if and only if $F^\perp$ is a flat vector subbundle of $(E,\na)$ and
the orthogonal decomposition $$E = F\oplus F^\perp$$ is a direct sum
of flat vector bundles.
\end{prop}

\section{Integration by parts} \label{int-by-parts-sec}
The main difference we will discuss between complex and affine
manifolds is in integration theory.  On an $n$-dimensional complex
manifold, an $(n,n)$ form is a volume form which can be integrated,
while on an affine manifold, an $(n,n)$ form is not a volume form.
Here we make a crucial extra assumption to handle this case
adequately:  We assume our affine manifold $M$ is equipped with a
$D$-invariant volume form $\nu$. Equivalently, we assume the linear
part of the holonomy of $D$ lies in $\Sl n$.  We call such an affine
manifold $(M,D,\nu)$ a \emph{special affine manifold}. This
important special case of affine manifold is quite commonly used: in
Strominger-Yau-Zaslow's conjecture \cite{syz}, a Calabi-Yau manifold
$N$ near the large complex structure limit in moduli should be the
total space of a (possibly singular) fibration with fibers of
special Lagrangian tori over a base manifold which is special
affine.  (The $D$-invariant volume form is the restriction of the
holomorphic $(n,0)$ form on $N$.)  Also, a famous conjecture of
Markus states that a compact affine manifold $(M,D)$ admits a
$D$-invariant volume form if and only if $D$ is complete.

From now on, we assume that $M$ admits a $D$-invariant volume
form $\nu$. Then $\nu$ provides natural maps from
\begin{eqnarray*} \mathcal A^{n,p} \to \Lambda^p, &&
\nu\otimes\chi\mapsto (-1)^{\frac{n(n-1)}2} \chi; \\ \mathcal
A^{p,n} \to \Lambda^p, && \chi\otimes\nu \mapsto
(-1)^{\frac{n(n-1)}2}\chi.
\end{eqnarray*}
(The choice of sign is to ensure that for every Riemannian metric
$g$, $\omega_g^n/\nu$ has the same orientation as $\nu$.) We use
division by $\nu$ to denote both of these maps. In particular,
$\chi\in\mathcal A^{n,n}$ can be integrated on $M$ by considering
$$\int_M \frac \chi\nu.$$

The reason we require $\nu$ to be $D$-invariant is to allow the
usual integration by parts formulas for $(p,q)$ forms to work
on the affine manifold $M$.  The main result we need is the
following:

\begin{prop} \label{int-by-parts}
Suppose $(M,D)$ is an affine flat manifold equipped with a
$D$-invariant volume form $\nu$. Then if $\chi\in \mathcal
A^{n-1,n}$, $$\frac {\partial \chi}\nu = d \left( \frac\chi
{2\nu} \right).$$ Also, if $\chi \in \mathcal A^{n,n-1}$,
$$\frac{\bar\partial \chi}{\nu} = (-1)^n\,\,d \left( \frac\chi
{2\nu} \right).$$
\end{prop}
\begin{proof}
We may choose local affine coordinates $x^1,\dots,x^n$ on $M$
so that $\nu = dx^1\wedge \dots \wedge dx^n$, and write
$\chi\in \mathcal A^{n-1,n}$ locally as
\begin{eqnarray*}
 \chi &=&
\sum_{i=1}^n f_i\, dz^1\wedge \cdots \wedge \widehat{dz^i} \wedge
\cdots \wedge dz^n \otimes d\bar z^1 \wedge \cdots \wedge d\bar
z^n, \\
\partial \chi &=& \frac12 \sum_{i=1}^n (-1)^{i-1} \frac{\partial f_i}{\partial
x^i} \,dz^1\wedge  \cdots
\wedge dz^n \otimes d\bar z^1 \wedge \cdots \wedge d\bar z^n,\\
 \frac \chi \nu &=& (-1)^{\frac{n(n-1)}2} \sum_{i=1}^n f_i \,dx^1\wedge \cdots
\wedge \widehat{dx^i} \wedge \cdots \wedge dx^n, \\
d\left(\frac\chi\nu\right) &=& (-1)^{\frac{n(n-1)}2}
\sum_{i=1}^n(-1)^{i-1} \frac{\partial f_i}{\partial x^i}\, dx^1
\wedge \cdots \wedge dx^n.
\end{eqnarray*}
(Note that when restricted to $M$, $dz^i=d\bar z^i = dx^i$.)
The computation is similar for $\chi\in\mathcal A^{n,n-1}$.
\end{proof}

%\begin{rem}
%The flat connection $D$ does need to be torsion-free in order
%for this result to hold.  If $D$ is flat but not torsion-free,
%then we may still define  $\mathcal A^{p,q}$, $\partial$, and
%$\bar\partial$ similarly, but other terms come into the
%calculation.
%\end{rem}

To each  Riemannian metric $g$ on an affine manifold $M$, there is a
natural nondegenerate $(1,1)$ form given by $\omega_g= g_{ij}\,
dx^i\otimes dx^j$ for $x^i$ local coordinates on $M$.  (The metric
$g$ is naturally extended to a Hermitian metric on $M^\co$ and
$\omega_g$ is the restriction of the Hermitian form of $g$ to
$M\subset M^\co$.) A metric $g$ on $M$ is said to be \emph{affine
Gauduchon} if $\partial\bar\partial (\omega_g^{n-1}) = 0$.  We will
see in the next section that every conformal class of Riemannian
metrics on a compact special affine manifold $M$ contains an affine
Gauduchon metric.

Note that by our convention (\ref{form-convention}) our definition
of first Chern form is a \emph{real} $(1,1)$ form on $M$, even
though it is the restriction of an \emph{imaginary} 2 form on
$M^\co$.

A locally constant vector bundle $E$ over a special affine manifold
$(M,\nu)$ equipped with an affine Gauduchon metric $g$ has a
\emph{degree} given by
\begin{equation}\label{degree-def}
\deg_g E = \int_M \frac{c_1(E,h) \wedge \omega_g^{n-1}} \nu.
\end{equation}
Recall the affine first Chern form $c_1(E,h)=-\partial \bar\partial
\log \det h_{\alpha\bar\beta}$ for any Hermitian metric $h$ on $E$.
The degree is well-defined because
 \begin{itemize}
\item For any other metric $h'$ on $E$,
$$c_1(E,h')-c_1(E,h)= \partial\bar\partial
(\log \det h_{\alpha\bar \beta} - \log \det h'_{\alpha \bar
\beta}),$$ which is $\partial\bar\partial$ of a function on $M$.
 \item Proposition \ref{int-by-parts} above allows us to integrate
 by parts to move the $\partial\bar\partial$ to $\omega_g^{n-1}$.
 \item The metric $g$ is affine Gauduchon.
\end{itemize}

Note we do not expect the degree to be an integer (see e.g.\
L\"ubke-Teleman \cite{lubke-teleman95} for counterexamples in
the complex case).

The \emph{slope} of a flat vector bundle $E$ over a special affine
manifold $M$ equipped with an affine Gauduchon metric $g$ is defined
to be
$$ \mu_g = \frac{{\rm deg}_g E } {{\rm rank}\,E}.$$
Such a complex flat vector bundle $E$ is called \emph{$\co$-stable}
if every flat subbundle $F$ of $E$ satisfies
 \begin{equation} \label{slope-stable} \mu_g(F) < \mu_g(E).
\end{equation} A real flat vector bundle $E$ is called
\emph{$\re$-stable} if (\ref{slope-stable}) is satisfied by any
flat real vector subbundle $F$ of $E$, while such an $E$ is
called $\co$-stable if the complex vector bundle
$E\otimes_\re\co$ is $\co$-stable.

\section{Affine Hermitian-Einstein metrics} \label{ahem-sec}

In this section, we will check some of the basic properties of
Hermitian-Einstein metrics extend to the affine case: a
vanishing theorem of Kobayashi, uniqueness of affine
Hermitian-Einstein metrics on simple bundles, and stability of
flat bundles admitting affine Hermitian-Einstein metrics. These
roughly form the easy part of the Kobayashi-Hitchin
correspondence between stable bundles and Hermitian-Einstein
metrics.  The hard part, to prove the existence of
Hermitian-Einstein metrics, will be addressed in the Sections
\ref{cont-meth-sec} to \ref{destabilizing-sec} below.

We have the following vanishing theorem of Kobayashi
\cite{kobayashi87}
\begin{thm} \label{van-thm}
Let $(E,\na)$ be a flat vector bundle over a compact affine manifold
$M$ equipped with a Riemannian metric $g$.  Assume $E$ admits an
affine Hermitian-Einstein metric $h$ with Einstein factor
$\gamma_h$.  If $\gamma_h<0$, then $E$ has no nontrivial locally
constant sections. If $\gamma_h=0$, then any locally constant
section $s$ of $E$ satisfies $\partial^h s=0$ for $\partial^h =
\partial^{\na,h}$.
\end{thm}

\begin{proof}
For $s$ any locally constant section of $E$, compute $${\rm tr}_g
\partial \bar
\partial |s|^2 = - \gamma_h |s|^2 + |\partial^h s|^2$$ and apply the
maximum principle.
\end{proof}

The following uniqueness proposition follows L\"ubke-Teleman
\cite[Prop.\ 2.2.2]{lubke-teleman95}
\begin{prop}
If $(E,\na)$ is a simple flat vector bundle over a compact affine
manifold $M$ with Riemannian metric $g$, then any
$g$-affine-Hermitian-Einstein metric on $E$ is unique up to a
positive scalar.
\end{prop}

\begin{proof}
Let $h_1,h_2$ be two affine Hermitian-Einstein metrics on $E$ with
Einstein constants $\gamma_1,\gamma_2$.  Then there an endomorphism
$f$ of $E$ satisfying $h_2(s,t)=h_1(f(s),t)$ for all sections $s,t$,
and since $h_1,h_2$ are both Hermitian, $f^\frac12$ is well-defined.

Then the connection $\na'=f^\frac12 \circ \na \circ f^{-\frac12}$ is
a flat connection on $E$.  Let $E'$ signify the new flat structure
$\na'$ induces on the underlying vector bundle of $E$, and let $E$
signify the original flat structure $\na$.  Then $f^\frac12$ is a
locally constant section of the flat vector bundle ${\rm
Hom}(E,E')$, $h_1$ is affine Hermitian-Einstein with Einstein
constant $\gamma_2$ on $E'$, and so the metric induced on ${\rm
Hom}(E,E')$ by $h_1$ on $E'$ and $h_2$ on $E$ is affine
Hermitian-Einstein with Einstein constant $\gamma_2-\gamma_2=0$.

Therefore, Theorem \ref{van-thm} applies, to show that
$\partial_{\rm Hom} (f^\frac12) = 0$ for $\partial_{\rm Hom}$ the
$(1,0)$ part of the affine Hermitian-Einstein connection on ${\rm
Hom}(E,E')$.  A computation as in \cite{lubke-teleman95} then
implies that $\partial_1 f = 0$ for $\partial_1$ the $(1,0)$ part of
the affine Hermitian connection on $(E, h_1)$. Since $f$ is
$h_1$-self-adjoint, this implies $\bar\partial (f^*) = \bar \partial
f = 0$.

So since $(E,\na)$ is simple, $f$ is a multiple of the identity.
\end{proof}

The following theorem is due to Kobayashi in the K\"ahler case
\cite{kobayashi87}. The proof in the present case is simpler
because we need only deal with subbundles and not singular
subsheaves in the definition of stability.
\begin{thm}
Let $E$ be a flat vector bundle over a compact special affine
manifold $M$ equipped with an affine Gauduchon metric $g$.  If $E$
admits an affine Hermitian-Einstein metric $h$ with Einstein
constant $\gamma$, then either $E$ is $g$-stable or $E$ is an
$h$-orthogonal direct sum of flat stable vector bundles, each of
which is affine Hermitian-Einstein with Einstein constant $\gamma$.
\end{thm}

\begin{proof}
Consider $E'$ a flat subbundle of $E$.  Then it suffices to prove
that $\mu(E')\le \mu(E)$ with equality only in the case that the
$h$-orthogonal complement of $E'\subset E$ is also a flat subbundle
of $E$.  By Proposition \ref{second-fund-prop} above, it suffices to
show that $\mu(E')\le \mu(E)$ with equality only if the second
fundamental form of $E'\subset E$ vanishes.

We compute, as in \cite[Proposition 2.3.1]{lubke-teleman95} or
\cite[Proposition V.8.2]{kobayashi87} that for $s= {\rm rank}\, E'$,
$r = {\rm rank}\,E$, that
\begin{eqnarray*}
 \mu_g E &=& \frac1{rn} \int_M {\rm tr}\, K_E \frac{\omega_g^n}\nu \\
 &=& \frac\gamma n\int_M \frac{\omega_g^n}\nu , \\
 \mu_g E' &=& \frac1{sn}\int_M {\rm tr}\, K_{E'} \frac{\omega_g^n}
 \nu \\
 &=& \frac{\gamma} n \int_M \frac{\omega_g^n}\nu - \frac1{sn} \int_M
 |A|^2 \frac{\omega_g^n}\nu.
\end{eqnarray*}
Thus $\mu_g E\le \mu_g E'$ always, with equality if and only if the
second fundamental form $A$ is identically 0.

For the exact sequence of flat bundles $0\to E' \to E \to E''\to 0$,
the extended curvatures $R'$, $R$, and $R''$ of the Hermitian
connections induced by $h$ on $E'$, $E$, $E''$ respectively, satisfy
$$ R  =\left( \begin{array}{cc} R' + A \wedge A^* & * \\
* & R'' + A^* \wedge A \end{array} \right)$$ (see e.g.\ Kobayashi
\cite[Proposition I.6.14]{kobayashi87}).  So the vanishing of $A$
implies that the mean curvatures $K' = {\rm tr}_g R'$ of $E'$ and
$K'' = {\rm tr}_g R''$ of $E''$ satisfy the Hermitian-Einstein
condition if $K$ does.  Thus, in the case of equality $\mu_g E =
\mu_g E'$, $E$ splits into proper flat affine Hermitian-Einstein
summands $E'$ and $E''$. The theorem then follows by induction on
the rank $r$.
\end{proof}

\section{Affine Gauduchon metrics}
\label{aff-g-section}

Given a smooth Riemannian metric $g$ on an affine manifold $M$ with
parallel volume form $\nu$, define the operator from functions to
functions given by
\begin{equation} \label{Q-def}
Q(\phi) = \frac{\partial\bar\partial(\phi\, \omega_g^{n-1})}
{\omega_g^n}.
\end{equation}
If we can find a smooth, positive solution to $Q(\phi)=0$, then
$\phi^{\frac1{n-1}} g$ is affine Gauduchon.

Consider the adjoint $Q^*$ of $Q$ with respect to the inner product
\begin{equation}\label{in-prod}
\langle \phi,\psi \rangle_g = \int_M \phi\,\psi \frac{\omega_g^n}
\nu.
 \end{equation}
Note that we are \emph{not} integrating with respect to the
volume form of $g$.  We can avoid extra curvature terms by
using the volume form $\omega_g^n/\nu$ instead (these terms are
worked out in the case of affine K\"ahler manifolds by Shima
\cite{shima07}). Compute, using Proposition \ref{int-by-parts}
above,
\begin{eqnarray*}
\langle \phi, Q^*(\psi)\rangle_g &=& \langle Q(\phi), \psi \rangle_g \\
&=& \int_M \frac{\partial \bar \partial (\phi\,\omega_g^{n-1})}
{\omega_g^n} \, \psi \, \frac{\omega_g^n}{\nu} \\
&=& \int_M \phi\, \frac{\partial\bar \partial \psi \wedge
\omega_g^{n-1}} \nu, \\
Q^*(\psi) &=& \frac{\partial\bar\partial \psi\wedge \omega_g^{n-1}}
{\omega_g^n} \\
&=& \frac1{4n} \, g^{ij} \frac{\partial^2\psi}{\partial x^i
\partial x^j}
= \frac1n \mbox{tr}_g \partial \bar \partial \psi.
\end{eqnarray*}
We have the following lemma
\begin{lem} \label{Q-star-max-prin}
The kernel of $Q^*$ consists of only the constant functions.  The
only nonnegative function in the image of $Q^*$ is the zero
function.
\end{lem}
\begin{proof}
Both statements follow directly from the strong maximum principle.
\end{proof}

The index of $Q$ (and of $Q^*$) is 0, as it is an elliptic
second-order operator on functions. The previous lemma shows
the kernel of $Q^*$ is one-dimensional, and thus the cokernel
of $Q^*$ (which may be identified with the kernel of $Q$ by
orthogonal projection) is one-dimensional as well. We want to
exhibit a positive function in the one-dimensional space $\ker
Q$.

Let $\phi\in \ker Q$ be not identically zero. If $\psi$ is not
in the image of $Q^*$, then $\langle \phi,\psi \rangle_g \neq
0$. This is because the dimension of the cokernel of $Q^*$ is
one, and the functional
$$\psi \mapsto \langle \phi,\psi\rangle_g$$ is not identically zero
but is zero on the image of $Q^*$. If $\phi$ assumes both positive
and negative values, then we can find a positive function $\psi$ on
$M$ so that $\langle \phi,\psi \rangle_g = 0$. But Lemma
\ref{Q-star-max-prin} above shows this $\psi$ is not in the range of
$Q^*$, a contradiction. Therefore, $\phi$ does not assume both
positive and negative values. Assume without loss of generality that
$\phi\ge0$.

Now, since $\phi\in \ker Q$ is not identically zero, and since $Q$
is an elliptic linear operator, the strong maximum principle shows
that $\phi>0$.  $C^\infty$ regularity of $\phi$ is standard. So the
above discussion has proved

\begin{thm}
If $M$ is a compact affine manifold with covariant-constant volume
form $\nu$, then every conformal class of Riemannian metrics on $M$
contains an affine Gauduchon metric unique up to scaling by a
constant.
\end{thm}

We will need the following lemma later.
\begin{lem} \label{Q-kernel-const}
If $g$ is an affine Gauduchon metric on a compact special
affine manifold, then the kernel of $Q$ consists only of the
constant functions.
\end{lem}
\begin{proof}
 If $\partial \bar \partial \omega_g^{n-1} = 0$, then the definition
 (\ref{Q-def}) shows that in local affine coordinates, $Q$
 is an elliptic operator of the form
 $$Q(\phi) = a^{ij} \frac{\partial^2\phi}{\partial x^i \partial x^j}
 + b^j \frac{\partial \phi}{\partial x^j}.$$
 So the strong maximum principle applies, and any function in the
 kernel of $Q$ must be constant.
\end{proof}

\section{The continuity method} \label{cont-meth-sec}

Consider a compact affine manifold $M$ equipped with a
covariant-constant volume form $\nu$ and an affine Gauduchon metric
$g$, and a flat complex vector bundle $E$ over $M$, together with a
Hermitian metric $h_0$.  Let $K_0$ be the extended mean curvature of
$(E,h_0)$.  Equations (\ref{trace-K}) and (\ref{degree-def}) show
$$ \int_M ({\rm tr}\, K_0)\, \frac{\omega_g^n}\nu = n\,\deg_g E,$$
and therefore for any affine Hermitian-Einstein metric on $E$
(satisfying $K=\gamma\,I_E$), $\gamma$ must satisfy
 \begin{equation}\label{gamma-def}
\gamma  \int_M \frac{\omega_g^n} \nu   = n\,\frac{\deg_g E}
{\mbox{rank}\, E} =  n\, \mu_g E .
\end{equation}

Let $h_0$ be a background Hermitian metric  $E$. Then any other
Hermitian metric $h$ on $E$ is given may be represented by an
endomorphism $f$ of $E$, so that for sections $s,t$, $$ h(s,t) =
h_0(f(s),t) \qquad \Longleftrightarrow \qquad f_\alpha^\eta =
h_0^{\eta \bar \beta} h_{\alpha\bar \beta}.$$ The new metric $h$ is
Hermitian if and only if $f$ is Hermitian self-adjoint and positive
with respect to $h_0$.  Here are some standard formulas for how the
extended connection form $\theta$, curvature $\Omega$,  first Chern
form $c_1$ and mean curvature $K$ change when passing from $h_0$ to
$h$:
\begin{eqnarray}
\theta &=& \theta_0 + f^{-1} \partial_0 f, \\
\Omega &=& \bar\partial\theta = \Omega_0 + \bar\partial (f^{-1}
\partial_0 f), \\
 \label{change-K}
K &=& K_0 + \mbox{tr}_g [\bar\partial (f^{-1}
\partial_0 f)], \\
c_1(E,h) &=& c_1(E,h_0) - \partial\bar \partial \log\det f, \\
\label{change-tr-K}
 \mbox{tr}\, K &=& \mbox{tr}\, K_0 -
\mbox{tr}_g \partial\bar\partial(\log \det f).
\end{eqnarray}
Note that in a locally constant frame, $f^{-1}\partial_0 f$ may be
written as $(f^{-1})^\alpha_\eta (\partial_0 f)^\eta_\beta$. The
term $\partial_0f$ is the extended Hermitian connection induced from
$(E,h_0)$ onto $\mbox{End}\,E$ acting on $f$:
$$(\partial_0 f)^\alpha_\beta = \partial f^\alpha_\beta -
(\theta_0)^\eta_\beta f^\alpha_\eta + (\theta_0)^\alpha _\eta
f^\eta_\beta.$$

Equation (\ref{change-K}) shows that we want to solve the equation
$$K_0 - \gamma \, I_E + \mbox{tr}_g[\bar \partial (f^{-1} \partial_0
f)] = 0.$$ We will solve this by the continuity method.  In
particular, for $\epsilon\in[0,1]$, consider the equation
\begin{equation} \label{cont-method-def}
L_\epsilon(f) = K_0 - \gamma \, I_E + \mbox{tr}_g[\bar \partial
(f^{-1} \partial_0 f)] + \epsilon\, \log f = 0.
\end{equation}
Note that since $f$ is an endomorphism of $E$ which is positive
Hermitian with respect to $h_0$, $\log f$ is well-defined.

Assume the background data $g$ and $h_0$ are smooth. Let $$J =
\{\epsilon\in (0,1] : \mbox{there is a smooth solution to }
L_\epsilon(f) = 0\}.$$ We will use the continuity method to show
that $J=(0,1]$ for any $\co$-simple flat vector bundle $E$, and then
later show that we may take $\epsilon\to0$ if $E$ is $\co$-stable.
(If $E$ is $\co$-stable, it is automatically $\co$-simple---see
Proposition \ref{co-simple} below.)

The first step in the continuity method is to show $1\in J$ and so
$J$ is nonempty.  The proof will also provide an appropriately
normalized initial metric $h_0$ on $E$.
\begin{prop} \label{est-init-cond}
Given a compact special affine manifold $M$ with an affine Gauduchon
metric and a flat vector bundle $E$. Then there is a smooth
Hermitian metric $h_0$ on $E$ so that there is a smooth solution
$f_1$ to $L_1(f)=0$.  The metric $h_0$ satisfies the normalization
${\mbox{tr}\, K_0} = r\,\gamma$ for $r$ the rank of $E$ and $\gamma$
given by (\ref{gamma-def}).
\end{prop}

\begin{proof}
We first produce the metric $h_1$ the metric satisfying the $L_1$
equation, and then we will produce $h_0$ from $h_1$.

Given an arbitrary background metric $h_0'$, equation
(\ref{change-tr-K}) above shows that if $h_1=e^\rho h_0'$ satisfies
${\mbox{tr}\, K_1} = r\,\gamma$ if and only if
 \begin{equation} \label{tr-K-const-eq}
 \mbox{tr}_g\partial\bar \partial \rho = \frac1{r} \,\mbox{tr}\, K_0' - \gamma
 \end{equation}
for $r$ the rank of $E$. Note the right-hand side satisfies
\begin{equation} \label{integral-0}
\int_M \left(\frac1{r} \,\mbox{tr}\, K_0' - \gamma \right)
\frac{\omega_g^n}\nu = 0.
 \end{equation}
Lemma \ref{Q-kernel-const} shows that the kernel of $Q$ consists of
only constants. Equation (\ref{integral-0}) then shows that the
right-hand side of (\ref{tr-K-const-eq}) is orthogonal to $\ker Q$
with respect to the inner product (\ref{in-prod}), and so must be in
the image of $Q^* = \frac1n \mbox{tr}_g \partial \bar \partial $.

Now define $f_1 = \exp(-K_1 + \gamma I_E)$ and
$$(h_0)_{\alpha\bar \beta} = (f_1^{-1})^\eta_\alpha
(h_1)_{\eta\bar \beta}.$$ Then we may check as in
L\"ubke-Teleman \cite[Lemma 3.2.1]{lubke-teleman95} that $h_0$
is a Hermitian metric and that, with respect to $h_0$,  $f_1$
satisfies $L_1(f_1)=0$. Moreover,
\begin{eqnarray*}
{\mbox{tr}\, K_0} &=& \mbox{tr}\, K_1 + \mbox{tr}_g
\partial \bar \partial \log \det f_1 \\
&=& \mbox{tr}\, K_1 + \mbox{tr}_g
\partial \bar \partial \,\mbox {tr} (-K_1 + \gamma I_E) \\
 &=& \mbox{tr}\, K_1 = r\,\gamma.
\end{eqnarray*}
\end{proof}

So for the choice of $h_0$ derived in Proposition
\ref{est-init-cond}, we have
\begin{cor} \label{nonempty-cor}
$1\in J$.
\end{cor}

\section{Openness}
Consider the $\mbox{Herm}(E,h_0)$ to be the space of endomorphisms
of the vector bundle $E$ which are Hermitian self-adjoint with
respect to $h_0$. In particular, we may check as in e.g.\
\cite[Lemma 3.2.3]{lubke-teleman95} that for $f$ a positive
Hermitian endomorphism of $E$, the operator
$$\hat L(\epsilon, f) = f L_\epsilon(f) = f K - \gamma f +  \epsilon f \log f
\in \mbox{Herm}(E,h_0).$$  Let $1<p<\infty$ and $k$ be a
sufficiently large integer.

Assume $\epsilon \in J$---in other words, there is a smooth solution
$f_\epsilon$ to $L_\epsilon(f) = 0 \Longleftrightarrow \hat
L(\epsilon, f) = 0$.  Then we will use the Implicit Function Theorem
to show that there is a $\delta>0$ so that for every $\epsilon'\in(
\epsilon - \delta, \epsilon + \delta)$, there is a solution to $\hat
L(\epsilon',f)=0$ in $L^p_k \mbox{Herm}(E,h_0)$.  Then, for $k$
large enough, we can bootstrap to show $C^\infty$ regularity of each
solution $f_{\epsilon'}$ to $\hat L(\epsilon', f_{\epsilon'}) = 0$.
Thus $( \epsilon - \delta, \epsilon + \delta) \cap (0,1] \subset J$
and $J$ is open.

So as usual, everything boils down to the checking the hypothesis of
the Implicit Function Theorem:
$$ \Xi = \frac{\delta}{\delta f} \,\hat L(\epsilon, f) \!: L^p_k
\mbox{Herm}(E,h_0) \to L^p_{k-2} \mbox{Herm}(E,h_0)$$ should be an
isomorphism of Banach spaces. The operator $\frac{\delta}{\delta f}
\,\hat L(\epsilon, f)$ is Fredholm and elliptic. The next thing to
check is that the index of the operator $\Xi$ is 0.
\begin{lem}
The index of $\Xi$ is 0.
\end{lem}
\begin{proof}
 To check this,
we need only look at the symbol.  For $\phi\in \mbox{Herm}(E,h_0)$,
compute
$$\Xi(\phi) \equiv \mbox{tr}_g \bar \partial \partial_0 \phi,$$ where $\equiv$
denotes equivalence up to zeroth- and first-order derivatives of
$\phi$.  Moreover, if $\phi,\xi\in \mbox{Herm}(E,h_0)$, then we may
compute
\begin{equation} \label{distr-deriv}
\bar\partial \left[ h_0(\partial_0 \phi,\xi) \right] =
h_0(\bar\partial
\partial_0 \phi, \xi) - h_0(\partial_0 \phi, \partial_0 \xi) .
\end{equation}
Here $h_0$ acts only on the $\mbox{End}(E)$ part of the quantities,
and not on the differential form parts: For $\phi_1,\phi_2$ sections
of $\mbox{End}(E)$, and $\lambda_i\in\mathcal A^{p_i,q_i}$,
\begin{equation} \label{h-forms}
h_0(\phi_1\otimes \lambda_1, \phi_2 \otimes \lambda_2) = h_0
(\phi_1,\phi_2) \lambda_1 \wedge \bar\lambda_2.
\end{equation}
The $\partial_0$ in the last time is because of the convention
(\ref{conj-form}) and the fact that $h_0$ is $\co$-antilinear in the
second slot, while the minus sign in front of the last term is
because of (\ref{h-forms}).

Now we use (\ref{distr-deriv}) to compute the highest-order terms of
the adjoint $\Xi^*$ of $\Xi$ with respect to the inner product
$$\langle \phi, \xi \rangle_{{\rm End}(E)} = \int_M
h_0(\phi,\xi)\frac{\omega_g^n} \nu.$$ Then compute using
(\ref{distr-deriv}) and Proposition \ref{int-by-parts}:
\begin{eqnarray*}
\langle \phi, \Xi^*\xi \rangle _{{\rm End}(E)} &=&  \langle
\Xi\phi,\xi \rangle_{{\rm End}(E)} \\
&=& \int_M h_0(\mbox{tr}_g\bar \partial \partial_0 \phi, \xi)
\frac{\omega_g^n} \nu \\
&=& n \int_M \frac{h_0(\bar \partial \partial_0 \phi, \xi) \wedge
\omega_g^{n-1}} \nu\\
&=& -n\int_M \frac{ h_0(\partial_0 \phi, \partial_0 \xi) \wedge
\omega_g^{n-1} - h_0(\partial_0\phi, \xi) \wedge \bar \partial
\omega_g^{n-1} } \nu \\
&=& n \int_M \frac{ h_0(\phi,\bar \partial \partial_0 \xi) \wedge
\omega_g ^{n-1} +T } \nu \\
&=&  \int_M \frac{h_0(\phi, \mbox{tr}_g \bar \partial \partial_0
\xi)\, \omega_g^n + T}\nu,
\end{eqnarray*}
where $T$ represents terms that involve no derivatives of $\phi$ and
only zeroth- or first-order derivatives of $\xi$.  Therefore, we see
$$\Xi^*(\phi) \equiv \mbox{tr}_g \bar\partial \partial_0 \phi
\equiv \Xi(\phi).$$ Since $\Xi$ and $\Xi^*$ have the same symbols,
they are homotopic as elliptic operators, and thus have the same
index.  Since the sum of the indices of $\Xi$ and $\Xi^*$ is 0, they
each must have index 0.
\end{proof}

Since the index of $\Xi$ is 0, it suffices to show $\Xi$ is
injective to apply the Implicit Function Theorem.  In order to do
this, we apply the following crucial estimate, essentially due to
Uhlenbeck-Yau.
\begin{prop} \label{eta-estimate}
Let $\alpha\in\re$, $\epsilon\in(0,1]$, $f$ be a positive and
Hermitian endomorphism of $E$ with respect to $h_0$, and $\phi\in
\mbox{Herm}(E,h_0)$. Assume $\hat L(\epsilon, f) = 0$ and $$
\frac{\delta}{\delta f} \hat L(\epsilon, f)(\phi) + \alpha f \log f
= \Xi(\phi) + \alpha f \log f =  0.$$ Then if $\eta = f^{-\frac12}
\phi f^{-\frac12},$
$$ -\mbox{tr}_g \partial\bar \partial |\eta|^2 + 2 \epsilon |\eta|^2
+ |\partial_0^f \eta|^2 + |\bar\partial^f \eta|^2 \le - 2 \alpha h_0
(\log f, \eta),$$ where $\partial_0^f = \mbox{Ad}\,f^{-\frac12}
\circ \partial_0 \circ \mbox{Ad}\, f^{\frac12}$ and $\bar \partial
^f = \mbox{Ad}\,f^{\frac12} \circ \bar\partial \circ \mbox{Ad}\,
f^{-\frac12}$, $|\partial_0^f\eta|^2 =\mbox{tr}_g h_0(\partial_0^f
\eta,\partial_0^f)$, and $|\bar\partial^f \eta|^2 = -\mbox{tr}_g
h_0(\bar \partial^f \eta, \bar \partial^f \eta)$.
\end{prop}
\begin{proof}
 This is a local calculation on $M$, which by our definitions of
 extended Hermitian connections, $p,q$ forms, etc., is the same as the
 calculation on $M^\co$.  So we refer the reader to \cite[Proposition 3.2.5]{lubke-teleman95}.
\end{proof}

\begin{prop} \label{open-prop}
$J$ is open.
\end{prop}
\begin{proof}
By the discussion above, we need only check that $\Xi$ in injective.
This follows from the previous Proposition \ref{eta-estimate} with
$\alpha=0$.  In this case, $$-\mbox{tr}_g \partial\bar \partial
|\eta|^2 + 2\epsilon |\eta|^2 \le 0,$$ and the maximum principle
implies $|\eta|^2=0$. So $\eta=0$ and $\phi=0$. $\Xi$ is injective.
\end{proof}

\section{Closedness}
\begin{lem} \label{det-1}
 If $f$ is a Hermitian positive endomorphism of $E$ with respect to $h_0$
 which solve $L_\epsilon(f)=0$ for $\epsilon>0$, then $\det f = 1$.
\end{lem}
\begin{proof}
Taking the trace of the definition (\ref{cont-method-def}) and using
Proposition \ref{est-init-cond}, we see that
$$-\mbox{tr}_g\partial \bar \partial \log \det f + \epsilon \log
\det f = 0.$$ The maximum principle then implies $\log \det f = 0$.
\end{proof}

We introduce some more notation.  Let $f=f_\epsilon$ represent the
family of solutions constructed for $\epsilon$ in the interval
$(\epsilon_0,1]$ in Corollary \ref{nonempty-cor} and Proposition
\ref{open-prop}. Define
$$ m= m_\epsilon = \max |\log f_\epsilon|, \quad \phi =
\phi_\epsilon = \frac{df_\epsilon} {d\epsilon}, \quad \eta =
\eta_\epsilon = f_\epsilon^{-\frac12} \phi_\epsilon f_\epsilon
^{-\frac12}.$$

We can immediately verify
\begin{lem}
The trace ${\rm tr}\,\eta_\epsilon=0$.
\end{lem}
\begin{proof}
Compute $$\mbox{tr}\,\eta = \mbox{tr}\, (f^{-\frac12}\phi
f^{-\frac12}) = \mbox{tr}\left(f^{-1}\frac{df}{d\epsilon} \right) =
\frac d{d\epsilon} (\log \det f) = 0
$$
by Lemma \ref{det-1} above.
\end{proof}

\begin{prop} \label{orthog-kernel}
Let $E$ be a $\co$-simple complex flat vector bundle over a compact
special affine manifold $M$.  On $M$, consider the $L^2$ inner
products on $\mathcal A^{p,q}(\mbox{End}\,E)$ given by $h_0$, $g$
and the volume form $\omega_g^n/\nu$. Then there is a constant
$C(m)$ depending only on $M$, $g$, $h_0$, $\nu$ and $m=m_\epsilon$
so that for $\eta=\eta_\epsilon$,
$$\|\bar\partial^f \eta \|^2_{L^2} \ge C(m) \| \eta \|^2_{L^2}.$$
\end{prop}

\begin{rem}
In the following sections, $C(m)$ will denote a constant depending
on $m$ and the other objects noted above, but the particular
constant may change with the context.  $C$ will similarly denote a
constant depending only on the initial conditions $M$, $g$, $h_0$
and $\nu$, but not on $\epsilon$ or $m$.
\end{rem}

\begin{proof}
 Let $\psi = f^{-\frac12}\eta f^{\frac12}$. Then pointwise, $$|\bar \partial^f
 \eta|^2 = |f^{\frac12} \bar \partial \psi f^{-\frac12}|^2 \ge
 C(m) |\bar \partial \psi|^2.$$ Integrate over $M$ with respect to
 the volume form $\omega_g^n/\nu$ to find that $$\| \bar \partial^f
 \eta \|^2_{L^2} \ge C(m) \| \bar \partial \psi \|^2 _{L^2} = C(m)
 \langle \bar \partial^* \bar \partial \psi, \psi \rangle,$$ where
 $\bar \partial^*$ is the adjoint of $\bar \partial$ with respect to
 the $L^2$ inner products on $\mathcal A^{0,0}(\mbox{End}\,E)$ and
 $\mathcal A^{0,1}(\mbox{End}\,E)$.  It is straightforward to check
 that $\bar\partial^*\bar\partial \!:\mathcal A^{0,0}(\mbox{End}
 \,E) \to \mathcal A^{0,0} (\mbox{End} \,E)$ is elliptic, and it is
 self-adjoint by formal properties of the adjoint.

 Now $\mbox{tr}\, \psi = \mbox{tr}\, (f^{-\frac12}\eta
 f^{\frac12}) = \mbox{tr}\, \eta = 0$, and so for $I_E$ the identity
 endomorphism of $E$,
 $$ \langle \psi, I_E \rangle_{L^2} = \int_M h_0(\psi,I_E) \frac{
 \omega_g^n} \nu = \int_M \mbox{tr} (\psi I_E) \frac{
 \omega_g^n} \nu = 0,$$
since $h_0(\psi,I_E) = \mbox{tr}(\psi I_E^*)$ for $I_E^*=I_E$ the
adjoint of $I_E$ with respect to $h_0$. Since $E$ is $\co$-simple,
this shows that $\psi$ is $L^2$-orthogonal to the kernel of
$\bar\partial$ on $\mbox{End}\,E$. Therefore, since $\bar\partial^*
\bar \partial$ is self-adjoint and elliptic, there is a constant
$\lambda_1>0$ (the smallest positive eigenvalue of $\bar\partial^*
\bar \partial$) so that
$$ \langle \bar\partial^* \bar \partial \psi, \psi \rangle_{L^2} \ge \lambda_1
\|\psi\|^2_{L^2}. $$ Therefore, $$\|\bar \partial^f \eta \|^2_{L^2}
\ge C(m)  \langle \bar\partial^* \bar \partial \psi, \psi
\rangle_{L^2} \ge C(m) \|\psi \|^2_{L^2} \ge C(m)
\|\eta\|^2_{L^2}.$$
\end{proof}

Now we need the following consequence of a subsolution estimate
of Trudinger \cite[Theorem 9.20]{gilbarg-trudinger}:
\begin{prop} \label{subsol-est}
If $u$ is a $C^2$ nonnegative function on $M$ which satisfies
$${\rm tr}_g \partial \bar \partial u \ge \lambda \,u + \mu$$ for
$\lambda\le 0$ and $\mu$ real constants, then
$$\max_M u \le B(\|u\|_{L^1} + |\mu|)$$ for $B$ a constant only
depending on $g$, $\nu$ and $\lambda$.
\end{prop}

Now we bound $|\phi|=|\phi_\epsilon|$ in terms of $m$.
\begin{prop} \label{sup-bound-phi}
Given $E$ a $\co$-simple complex flat vector bundle over a
compact special affine manifold $M$, $\max_M |\phi_\epsilon|
\le C(m)$.
\end{prop}
\begin{proof}
Proposition \ref{eta-estimate} above shows that $$-\mbox{tr}_g
\partial \bar \partial |\eta|^2 + |\bar\partial^f \eta|^2 \le 2
|\log f| \cdot |\eta|.$$ Since $\int_M\mbox{tr}_g\partial \bar
\partial |\eta|^2 \omega_g^n /\nu = 0$, we have $$\|\bar \partial ^f
\eta\|^2_{L^2} \le C(m) \,\|\eta\|_{L^2}.$$ But then Proposition
\ref{orthog-kernel} implies $$C(m)\,\|\eta\|^2_{L^2} \le \|\bar
\partial^f \eta \|^2_{L^2} \le C(m) \,\|\eta\|_{L^2} \qquad
\Longrightarrow \qquad \|\eta\|_{L^2}\le C(m).$$ But then we also
have from Proposition \ref{eta-estimate} that $$-\mbox{tr}_g
\partial \bar \partial |\eta|^2 \le 2 |\log f| \cdot |\eta| \le m \,
|\eta|^2 + m,$$ and Proposition \ref{subsol-est} then shows that
$$\max_M |\eta|^2 \le C(m)(\|\eta\|^2_{L^2} + m) \le C(m).$$
The result follows since $\phi = f^{\frac12}\eta f^{\frac12}$.
\end{proof}

The following lemma follows is a local calculation as in
\cite[Lemma 3.3.4.i]{lubke-teleman95}.
\begin{lem} \label{op-log-f}
 $$-\sfrac12 {\rm tr}_g \partial \bar \partial |\log f|^2 +
 \epsilon\, |\log f|^2 \le |K_0-\gamma I_E|\cdot |\log f|.$$
\end{lem}
\begin{cor} \label{m-ep-bound}
 $m \le \epsilon^{-1} C$.
\end{cor}
\begin{proof}
Apply the maximum principle to Lemma \ref{op-log-f} for $C=
\max_M|K_0- \gamma I_E|$.
\end{proof}
\begin{cor} \label{m-L2-bound}
$ m \le C (\|\log f\|_{L^2} + 1)^2.$
\end{cor}
\begin{proof}
Lemma \ref{op-log-f} implies $$ -\mbox{tr}_g \partial \bar \partial
|\log f|^2 \le |\log f|^2 + \max_M|K_0-\gamma I_E|^2.$$ Then
Proposition \ref{subsol-est} applies to show $$m \le C(\|\log
f\|^2_{L^2} + 1),$$ which implies the corollary.
\end{proof}

\begin{lem} \label{rel-to-self-adj}
Consider the operator $\bar\partial_0^* \bar \partial_0$ acting
on sections of $\mbox{End}(E)$, where the adjoint is with
respect to the inner product $\langle\cdot,\cdot\rangle_{{\rm
End}(E)}$. Then for each section $\psi$  of $\mbox{End}(E)$,
$$\partial_0^*
\partial_0 \psi = \frac1n\,{\rm tr}_g \bar\partial \partial_0 \psi -\,
\frac{\partial_0 \psi \wedge \bar \partial \omega_g^{n-1}}
{\omega_g^n}.$$
\end{lem}
\begin{proof}
Since \begin{eqnarray*}\bar \partial [ h_0 (\partial_0
\psi_1,\psi_2) \wedge \omega_g^n] &=& [h_0(\bar\partial \partial_0
\psi_1, \psi_2) - h_0 (\partial_0\psi_1, \partial_0 \psi_2)] \wedge
\omega_g^{n-1} \\&&{}- h_0 (\partial_0 \psi_1, \psi_2) \wedge \bar
\partial \omega_g^{n-1},
\end{eqnarray*}
Proposition \ref{int-by-parts} and Stokes' Theorem show that
\begin{eqnarray*}
\int_M  h_0(\partial_0^* \partial_0 \psi_1, \psi_2) \frac{
\omega_g^n} \nu &=&
 \int_M \frac{h_0 (\partial_0 \psi_1, \partial_0 \psi_2) \wedge
 \omega_g ^{n-1}} \nu \\
 &=& \int_M \frac{h_0(\bar \partial \partial_0 \psi_1, \psi_2) \wedge
 \omega_g^{n-1}} \nu \\
 &&{} - \int_M \frac{h_0( \partial_0 \psi_1, \psi_2) \wedge \bar
 \partial \omega_g^{n-1}} \nu \\
 &=& \frac1n \int_M h_0(\mbox{tr}_g\bar \partial \partial_0 \psi_1, \psi_2)
 \frac{
 \omega_g^{n}} \nu \\
 &&{} - \int_M \frac{h_0( \partial_0 \psi_1, \psi_2) \wedge \bar
 \partial \omega_g^{n-1}} \nu \\
\end{eqnarray*}
\end{proof}

\begin{prop} Assume $E$ is a $\co$-simple complex flat vector
bundle over $M$ a compact special affine manifold. Suppose
there is an $m\in\re$ so that $m_\epsilon\le m$ for all
$\epsilon\in (\epsilon_0, 1]$. Then for all $p>1$ and
$\epsilon\in (\epsilon_0,1]$,
$$ \|\phi_\epsilon\|_{L^p_2} \le C(m) (1+ \|f_{\epsilon}\|_{L^p_2}
),$$ where $C(m)$ may depend on $p$ as well as $m$ and the initial
data.
\end{prop}
\begin{proof}
The variation $\phi=\phi_\epsilon$ satisfies \begin{eqnarray*} 0 &=&
\frac{\delta}{\delta f} \hat L(\epsilon,f)(\phi) + f\log f \\
&=& \phi[K_0 -\gamma I_E + \epsilon\log f + \mbox{tr}_g \bar
\partial (f^{-1} \partial_0 f)] \\
&&{}- f \,\mbox{tr}_g\bar\partial(f^{-1} \phi f^{-1} \partial_0 f) +
f\,\mbox{tr}_g \bar \partial(f^{-1} \partial_0\phi)  \\
&&{} + f\log f + \epsilon f \left(\frac{\delta} {\delta f} \log
f\right)(\phi).
 \end{eqnarray*}
One computes then that
\begin{eqnarray*}
\mbox{tr}_g (\bar\partial \partial_0 \phi) &=&
   -\phi(K_0-\gamma I_E + \epsilon \log f) - \mbox{tr}_g(\bar
\partial f \wedge f^{-1} \phi f^{-1} \partial_0 f) \\
&&{}+ \mbox{tr}_g (\bar\partial f \wedge f^{-1}\partial_0 \phi)
 + \mbox{tr}_g(\bar \partial \phi \wedge f^{-1} \partial_0 f) \\
 &&{}- f\log f - \epsilon f \left(\frac{\delta} {\delta f}
\log f\right)(\phi)
\end{eqnarray*}
Then  Lemma \ref{rel-to-self-adj} above shows that for the operator
$\Lambda = n\, \partial_0^* \partial_0 + I_E$
\begin{equation} \label{Lambda-eq}
\begin{array}{c@{\,}l}
\Lambda\phi \,=&
 -\phi[K_0-(\gamma +1)I_E + \epsilon \log f] - \mbox{tr}_g(\bar
\partial f \wedge f^{-1} \phi f^{-1} \partial_0 f) \\[1mm]
 &{}+ \mbox{tr}_g (\bar\partial f \wedge f^{-1}\partial_0
\phi)
 + \mbox{tr}_g(\bar \partial \phi \wedge f^{-1} \partial_0 f) \\
 &\D{}- f\log f - \epsilon f \left(\frac{\delta} {\delta f}
\log f\right)(\phi) - n\,\frac{\partial_0 \phi \wedge \bar \partial
\omega_g^{n-1}} {\omega_g^n}.
\end{array}
\end{equation}
The operator $\Lambda\!: L^p_2(\mbox{End}\,E) \to L^p (\mbox{End}\,
E)$ is elliptic, self-adjoint, and is continuously invertible, since
$\partial_0^* \partial_0$ has nonnegative spectrum.  Therefore,
there is a $C$ satisfying
$$ \|\phi\|_{L^p_2} \le C \| \Lambda \phi\|_{L^p},$$
where as usual $C$ depends only on the initial data and $p$.

So we consider the $L^p$ norms of the 7 terms on the right-hand side
of (\ref{Lambda-eq}): The first term is bounded by $C(m)$ by
Proposition \ref{sup-bound-phi}, and the fifth is also bounded by
$C(m)$. Proposition \ref{sup-bound-phi} and H\"older's inequality
shows the second term is bounded by $C(m) \|f\|^2_{L^{2p}_1}$. The
third and fourth terms are both bounded by $C(m) \|f\|_{L^{2p}_1} \|
\phi \|_{L^{2p}_1}$. A local computation shows the sixth term is
bounded by $C(m)$, and the last term is clearly bounded by $C\| \phi
\|_{L^{2p}_1}$. So, altogether,
$$ \|\phi\|_{L^p_2} \le C(m) ( 1 + \|\phi\|_{L^{2p}_1} +
\|\phi\|_{L^{2p}_1} \|f\|_{L^{2p}_1} + \|f\|_{L^{2p}_1}^2 ).$$

An interpolation inequality of Aubin \cite[Theorem 3.69]{aubin98}
states that
$$ \|\psi\|_{L^{2p}_1} \le C \|\psi\|_{L^\infty}^{\frac12} \| \psi
\|_{L^p_2}^{\frac12} + \| \psi \|_{L^{2p}}.$$ Since both
$\|f\|_{L^\infty}, \|\phi\|_{L^\infty} \le C(m)$, a  simple
computation allows us to prove the proposition.
\end{proof}

\begin{cor} \label{unif-Lp2-bound}
Assume there is a smooth family of solutions $f_\epsilon$ to
$L_\epsilon (f_\epsilon)=0$, and that there is a uniform $m$ so that
$m_\epsilon\le m$ for all $\epsilon \in (\epsilon_0,1]$. Then for
all $\epsilon\in (\epsilon_0,1]$, $\|f_\epsilon \|_{L^p_2} \le
C(m)$, where $C(m)$ does not depend on $\epsilon$.
\end{cor}

\begin{proof}
Since $\phi_\epsilon = \frac d{d\epsilon} f_\epsilon$,
$$\frac d{d\epsilon} \|f_\epsilon \|_{L^p_2} \ge -\| \phi_\epsilon
\|_{L^p_2} \ge -C(m) (1+ \|f_\epsilon \|_{L^p_2}).$$ Then simply
integrate this ordinary differential inequality.
\end{proof}

\begin{prop} Assume $E$ is a $\co$-simple flat complex vector
bundle over $M$ a compact special affine manifold. Then
$J=(0,1]$. Moreover, if $\|f_\epsilon\|_{L^2}$ is bounded
independently of $\epsilon\in(0,1]$, then there exists a smooth
solution $f_0$ to the Hermitian-Einstein equation $L_0(f_0)=0$.
\end{prop}

\begin{proof}
The first statement will follow if we can show $J$ is closed.  In
particular, all we need to show is that if  $J=(\epsilon_0, 1]$ for
$\epsilon_0>0$, then there is a smooth solution $f_{\epsilon_0}$ to
$L_{\epsilon_0}(f_{\epsilon_0}) = 0$.  Corollaries \ref{m-ep-bound}
and \ref{unif-Lp2-bound} and then shows there is a constant $C$
satisfying $\|f_\epsilon\|_{L^p_2} \le C$ for all $\epsilon\in
(\epsilon_0, 1]$.  We will use this uniform estimate below to show
the existence of $f_{\epsilon_0}$.

Under the hypotheses of the second statement of the proposition, on
the other hand, Corollaries \ref{m-L2-bound} and
\ref{unif-Lp2-bound} together show that there is a $C$ so that for
all $\epsilon\in(0,1]$, $\|f_\epsilon\|_{L^p_2} \le C$.

Therefore, to prove the whole proposition, we may assume that for
$\epsilon_0\in[0,1)$, there is a constant $C$ and a smooth family of
solutions $f_{\epsilon}$ of $L_\epsilon(f_\epsilon)=0$ exists and
satisfies $\|f_\epsilon\|_{L^p_2} \le C$.  We will find a sequence
$\epsilon_i\to \epsilon_0^+$ so that $f_{\epsilon_0} = \lim
f_{\epsilon_i}$ is the solution we require.

Choose $p>n$. In this case, $L^p_1$ maps compactly into $C^0$, and
so $\log \!: L^p_1(\mbox{End}\,E) \to L^p_1(\mbox{End}\,E)$ is
continuous and the product of two functions in $L^p_1$ are also in
$L^p_1$. (See e.g.\ \cite{lubke-teleman95}.)

The uniform $L^p_2$ bound implies there is a sequence $\epsilon_i\to
\epsilon_0$ so that $f_{\epsilon_i}\to f_{\epsilon_0}$ converges
weakly in $L^p_2$, and strongly in $L^p_1$ and $C^0$.  Then compute,
in the sense of distributions, for $\alpha$ a smooth section of
$\mbox{End}(E)$,
\begin{eqnarray*}
\langle L_{\epsilon_0}(f_{\epsilon_0}),\alpha\rangle _{{\rm End}
(E)} &=& \langle L_{\epsilon_0}(f_{\epsilon_0}) -
L_{\epsilon_i}(f_{\epsilon_i}) \rangle_{{\rm End}(E)}\\
&=& \int_M h_0(\mbox{tr}_g[\bar \partial(f_{\epsilon_0}^{-1}
\partial_0 f_{\epsilon_0} - f_{\epsilon_i}^{-1} \partial_0
f_{\epsilon_i} )], \alpha) \,\frac{\omega_g^n}\nu \\
&&{} + \int_M h_0(\epsilon_0\log f_{\epsilon_0} - \epsilon_i \log
f_{\epsilon_i}, \alpha) \, \frac{\omega_g^n} \nu
\end{eqnarray*}
The second term goes to zero as $\epsilon_i\to\epsilon_0$ since
$f_{\epsilon_i}\to f_{\epsilon_0}$ in $C^0$.  Using Proposition
\ref{int-by-parts}, the first term can be written as
$$ \begin{array}{c}\D n\int_M  \frac{h_0( f_{\epsilon_0}^{-1}
\partial_0 f_{\epsilon_0} - f_{\epsilon_i}^{-1} \partial_0
f_{\epsilon_i}, \partial_0 \alpha) \wedge \omega_g^{n-1} } \nu
 \\[1mm]
\D {} + n \int_M \frac {h_0( f_{\epsilon_0}^{-1}
\partial_0 f_{\epsilon_0} - f_{\epsilon_i}^{-1} \partial_0
f_{\epsilon_i}, \alpha) \wedge \bar \partial \omega_g ^{n-1} } \nu .
\end{array}$$
Both these terms converge to 0 since $f^{-1}_{\epsilon_i} \partial_0
f_{\epsilon_i} \to f^{-1}_{\epsilon_0} \partial_0 f_{\epsilon_0}$ in
$L^p$.  Therefore, $L_{\epsilon_0} (f_{\epsilon_0})=0$ in the sense
of distributions.

Now we can compute in much the same way, for $f_{\epsilon_0}\in
L^p_2$, $\mbox{tr}_g\bar \partial \partial_0 f_{\epsilon_0} \in
L^p_1$. Therefore,  $f_{\epsilon_0} \in L^p_3$, and we can bootstrap
further to show that $f_{\epsilon_0}$ is smooth and is a classical
solution to $L_{\epsilon_0} (f_{\epsilon_0})=0$.
\end{proof}

\section{Construction of a destabilizing subbundle}
\label{dest-sub-sec} \label{destabilizing-sec} In this section, we
will construct a destabilizing flat subbundle of $E$ if
$\limsup_\epsilon \|f_\epsilon\| _{L^2} = \infty$. For a sequence
$\epsilon_i\to 0$, we will rescale by the reciprocal $\rho_i$ of the
largest eigenvalue of $f_{\epsilon_i}$. Then we will show that  the
limit
$$\lim_{\sigma\to0} \lim_{i\to \infty} (\rho_i f_{\epsilon_i})
^\sigma$$ exists and all of its eigenvalues are 0 or 1. A projection
to the destabilizing subbundle will be given by $I_E$ minus this
limit.

\begin{prop} \label{lap-f-sigma-bound}
If $\epsilon>0$, $0<\sigma\le1$, and $f$ satisfies
$L_\epsilon(f)=0$, then
$$ -\frac1\sigma \,{\rm tr}_g \partial \bar \partial ( {\rm tr}
f^{\sigma} ) + \epsilon \, h_0(\log f, f^\sigma) + |f^{-\frac\sigma
2}
\partial_0 (f^\sigma)|^2 \le -h_0 (K_0-\gamma I_E, f^\sigma).
$$
\end{prop}

\begin{proof}
This is a local computation, for which we refer to \cite[Lemma
3.4.4]{lubke-teleman95}.
\end{proof}

To rescale $f_\epsilon$ properly, consider the largest eigenvalue
$\lambda(\epsilon,x)$ of $\log f_\epsilon(x)$ for $x\in M$, and
define
$$ M_\epsilon = \max_{x\in M} \lambda(\epsilon,x), \qquad
\rho_\epsilon = e^{-M_\epsilon}.$$  Then since $\det f_\epsilon =
1$, $\rho_\epsilon\le1$ and we have the following straightforward
lemma:

\begin{lem} \label{basic-rho-lemma}
Assume $\limsup_{\epsilon\to0} \|f_\epsilon\|_{L^2} = \infty$. Then
\begin{enumerate}
\item $\rho_\epsilon f_\epsilon \le I_E$.
\item For each $x\in M$, there is an eigenvalue of $\rho_\epsilon
f_\epsilon$ less than or equal to $\rho_\epsilon$.
\item $\max_M \rho_\epsilon|f_\epsilon| \ge 1$.
\item There is a sequence $\epsilon_i\to0$ so that
$\rho_{\epsilon_i}\to0$.
\end{enumerate}
\end{lem}

\begin{prop} \label{construct-f-0-infty}
 There is a subsequence $\epsilon_i\to0$ so that
 $\rho_{\epsilon_i}\to 0$ and so that $f_i = \rho_{\epsilon_i}
 f_{\epsilon_i}$ satisfies
 \begin{enumerate}
 \item $f_i$ converges weakly in $L^2_1$ to an $f_\infty\neq 0$.
 \item As $\sigma\to0$, $f^{\sigma}_\infty$ converges
     weakly in $L^2_1$ to $f^0_\infty$.
 \end{enumerate}
\end{prop}

\begin{proof}
First of all, note that since each $f^\sigma_\epsilon$ is
positive-definite and self-adjoint with respect to $h_0$,
\begin{equation} \label{norm-trace}
|f^\sigma_{\epsilon}| \le \mbox{tr}\,f^\sigma_\epsilon \le \sqrt r\,
|f^\sigma_\epsilon|.
\end{equation}
Let $\sigma\in(0,1]$. Then Proposition \ref{lap-f-sigma-bound},
Corollary \ref{m-ep-bound}, and (\ref{norm-trace}) show
\begin{eqnarray*}
 \mbox{tr}_g\partial \bar \partial (\mbox{tr}\,f_\epsilon^\sigma) &\ge& \epsilon \,
 h_0(\log f_\epsilon,f_\epsilon^\sigma) + h_0(K_0-\gamma I_E,f_\epsilon^\sigma) \\
 &\ge& -(\epsilon m_\epsilon + C)|f_\epsilon^\sigma| \\
 &\ge& -C|f_\epsilon| \ge -C \,\mbox{tr}\,f_\epsilon^\sigma,
\end{eqnarray*}
where, as usual, $C$ is a (changing) constant depending only on the
initial data.  Now Proposition \ref{subsol-est}, Lemma
\ref{basic-rho-lemma} and (\ref{norm-trace}) show that
 \begin{equation} \label{lim-not-zero} 1 \le \max_M
\rho_\epsilon^\sigma |f_\epsilon^\sigma| \le \max_M
\rho_\epsilon^\sigma \mbox{tr}\,f_\epsilon^\sigma \le
C\rho_\epsilon^\sigma \|\mbox{tr}\,f_\epsilon^\sigma\|_{L^1} \le C
\|\rho_\epsilon^\sigma f_\epsilon^\sigma\|_{L^2}.
\end{equation}
On the other hand, Lemma \ref{basic-rho-lemma} shows $$\|\rho_
\epsilon^\sigma f_\epsilon^\sigma \|_{L^2} \le \|I_E\|_{L^2} = C,$$
and so it remains to estimate on $\|\partial_0 (f_i^\sigma)\|_{L^2}$
to get uniform bounds on $\|f_i^\sigma\|_{L^2_1}$.

Compute for $\epsilon=\epsilon_i$,
\begin{eqnarray*}
\|\partial_0 f_i^\sigma\|_{L^2}^2 &=& \int_M |\partial_0
(\rho_{\epsilon}^\sigma
f_{\epsilon}^\sigma)|^2 \, \frac{\omega_g^n}\nu \\
&\le& \int_M |(\rho_\epsilon f_\epsilon)^{-\frac\sigma2}
\partial_0
(\rho_{\epsilon}^\sigma f_{\epsilon}^\sigma)|^2 \, \frac{\omega_g^n}\nu \\
&\le& \rho_\epsilon^\sigma \int_M \frac1\sigma \mbox{tr}_g\partial
\bar
\partial (\mbox{tr} f^\sigma_\epsilon) \frac{\omega_g^n}\nu - \rho_\epsilon^\sigma \int_M
h_0(\epsilon\log
f_{\epsilon} + K_0 -\gamma I_E, f^\sigma_\epsilon) \, \frac{\omega_g^n}\nu \\
&=& \frac {\rho_\epsilon^\sigma n}\sigma \int_M \frac{\partial \bar
\partial (\mbox{tr} f^\sigma_\epsilon) \wedge \omega_g^{n-1}} \nu -
\int_M h_0(\epsilon\log f_{\epsilon} + K_0 -\gamma I_E,
\rho_\epsilon^\sigma f^\sigma_\epsilon)
\, \frac{\omega_g^n}\nu \\
&=& - \int_M h_0(\epsilon\log f_{\epsilon} + K_0 -\gamma I_E,
\rho_\epsilon^\sigma f^\sigma_\epsilon)
\, \frac{\omega_g^n}\nu \\
&\le & C \max_M (\rho_\epsilon f_\epsilon)^\sigma  \le C,
\end{eqnarray*}
where we have used Lemma \ref{basic-rho-lemma} to show
$(\rho_\epsilon f_\epsilon)^{-\frac\sigma2} \ge I_E$ to derive the
second line from the first; Proposition \ref{lap-f-sigma-bound} for
the third line; Proposition \ref{int-by-parts}, Stokes' Theorem, and
the fact that $g$ is affine Gauduchon to get the fifth line; and
finally Corollary \ref{m-ep-bound} and Lemma \ref{basic-rho-lemma}
to derive the sixth line. Note the final bound $C$ is independent of
$\sigma$ and $\epsilon$.

For $\sigma=1$, therefore, we have uniform $L^2_1$ bounds on $f_i$,
and so there is an $L^2_1$-weakly-convergent subsequence which we
may assume converges in $L^2$ and almost everywhere on $M$. For
simplicity, we still call this subsequence $f_i$. The bound
(\ref{lim-not-zero}) shows that $f_\infty = \lim f_i$ is not zero in
$L^2$.

The almost everywhere convergence of $f_i\to f_\infty$ shows that
$f_\infty$ is $h_0$-adjoint and positive semidefinite almost
everywhere.  Lemma \ref{basic-rho-lemma} shows that each eigenvalue
of $f_\infty$ is in $[0,1]$.  Therefore, by considering a
(measurable) frame which diagonalizes $f_\infty$ at almost every
point, it is clear that $f^\sigma_\infty$ converges to a limit
$f^0_\infty$ pointwise almost everywhere as $\sigma\to0$.

Moreover, the uniform bounds on $\|f^\sigma_i\|_{L^2_1}$ for
all $\sigma\in(0,1]$ show that $\|f_\infty^\sigma\|_{L^2_1}$ is
also uniformly bounded independent of $\sigma$, and so for each
sequence $\sigma_j\to0$, there is a subsequence $\sigma_{j_k}$
so that $f_\infty^{\sigma_{j_k}}$ converges weakly in $L^2_1$,
strongly in $L^2$ and pointwise almost everywhere to
$f^0_\infty$.  Thus $f^\sigma_\infty \to f^0_\infty$ weakly in
$L^2_1$ as $\sigma\to0$.
\end{proof}

Now let $\pi = I_E-f_\infty^0$.

\begin{prop} \label{produce-pi}
The endomorphism $\pi = I_E-f^0_\infty$ is an $h_0$-orthogonal
projection onto a flat subbundle of $E$.  In other words, it
satisfies $\pi^2=\pi$, $\pi^* = \pi$ and $(I_E-\pi)\bar
\partial \pi=0$ in $L^1$.  Moreover, $\pi$ is a smooth
endomorphism of $E$. So the locally constant subbundle
$F=\pi(E)$ is smooth.
\end{prop}

\begin{proof}
First we show that $\pi=\pi^*$, $\pi=\pi^2$, and $(1-\pi)\bar
\partial \pi=0$ in $L^1$ only. Then we will finish the proof with a discussion of
regularity.

To show $\pi=\pi^*$ almost everywhere, recall $f_\infty^0$ is a
pointwise almost-everywhere limit of $f_\infty^\sigma$, and
$f_\infty$ is a pointwise almost-everywhere limit of $f_i$, which
satisfies $f_i = f_i^*$.

To show $\pi^2=\pi$ in $L^1$, use Proposition
\ref{construct-f-0-infty} to compute
$$\pi^2 = \lim_{\sigma\to0} (I_E - f_\infty^\sigma)^2 = I_E -
2\lim_{\sigma\to0} (f_\infty^\sigma +  f_\infty^{2\sigma}) = 1 - 2
f^0_\infty + f^0_\infty = \pi.$$

To show $(1-\pi)\bar \partial \pi =0$ in $L^1$, compute since
$\pi=\pi^*=\pi^2$ that
$$ |(I_E-\pi)\bar\partial \pi| = |\bar \partial (I_E-\pi) \pi|
= |[\bar\partial(I_E-\pi) \pi]^* | = |\pi \partial_0(I_E - \pi)|.$$
(Here ${}^*$ represents the adjoint with respect to $h_0$ only, and
not with respect to any Hodge-type star on the affine Dolbeault
complex $\mathcal A^{p,q}({\rm End}\,E)$.) So we will show that
$$\|\pi\partial_0 (I_E-\pi)\|_{L^2} = 0.$$  Since the eigenvalues
of $f_i$ are between 0 and 1, a local computation (see e.g.\
\cite[p.\ 87]{lubke-teleman95}) implies that
$$0 \le \frac{s + \frac\sigma2}s (I_E - f_i^s) \le f_i
^{-\frac\sigma2}$$ for $0\le s\le \frac\sigma2$.  Then, as above,
Proposition \ref{lap-f-sigma-bound} shows that \begin{eqnarray*}
 \int_M | (I_E - f_i^s) \partial_0( f_i^\sigma)|^2 \frac{
 \omega_g^n} \nu &\le& \left(\frac s{s+\frac\sigma2}\right)^2
\int_M |f_i^{-\frac\sigma2} \partial_0 (f_i^\sigma)|^2 \frac{
 \omega_g^n} \nu \\\
 &\le&  \left(\frac s{s+\frac\sigma2}\right)^2
  \int_M |\epsilon_i \log f_i +K_0 - \gamma I_E|
 |f_i|^\sigma \frac{
 \omega_g^n} \nu\\
 &\le&\left(\frac s{s+\frac\sigma2} \right)^2  C.
\end{eqnarray*}
Since $\{(I_E-f_i^s) \partial_0(f^\sigma_i)\}^\infty_{i=1}$ is
a bounded sequence in $L^2$, weak compactness in $L^2$ allows
us to take $i\to\infty$ to find
$$\int_M | (I_E - f_\infty^s) \partial_0( f_\infty^\sigma)|^2 \frac{
 \omega_g^n} \nu \le \left(\frac s{s+\frac\sigma2} \right)^2  C.$$
Now we let $s\to0$ first so that $I_E-f_\infty^s \to
I_E-f_\infty^0=\pi$ strongly in $L^2$ as $s\to0$ by the uniform
$L^2_1$ bounds.   So $$ \int_M |\pi
\partial_0 (f^\sigma_\infty) | \frac{\omega_g^n} \nu = 0.$$
By definition, $\lim_{\sigma\to 0} \partial_0 f_\infty^\sigma$
converges weakly in $L^2$ to $\partial_0 (I_E-\pi)$, and so $\int_M
|\pi\partial_0(I_E-\pi)| \frac{\omega_g^n} \nu = 0$.

It remains to show that $\pi=\pi^2=\pi^*$ and $\pi \bar\partial
(I_E- \pi)=0$ in $L^1$ implies that $\pi$ is smooth.  The regularity
of $F=\pi(E)$ is a local issue, and so we restrict to a local
coordinate chart and a locally constant frame.  By an argument of
Popovici \cite[Lemma 0.3.3]{popovici05}, we can assume $h_0$ is the
standard flat metric with regards to the locally constant frame.

In terms of the standard flat metric, in order to show that
$F=\pi(E)$ is a smooth flat vector bundle, it suffices to show
that
$$\bar\partial \pi = 0 \qquad \Longleftrightarrow \qquad \nabla
\pi = 0.$$

At each $x\in M$, $\pi(x)$ can be considered as a map from $\co^r$
to $\co^r$ of some rank $k$.  The conditions $\pi$ satisfies are
then
$$ \pi^2 = \pi,\qquad \pi^* = \pi, \qquad (I_E-\pi)\bar \partial \pi
= 0,$$ for ${}^*$ the conjugate transpose.
 Now $\pi$ is $L^2_1$ when restricted to almost every coordinate
 line segment, with variable $t$ on the segment.  Then the last
 condition on $\pi$ becomes $$ (I-\pi)\frac{d\pi}
 {dt}= (I-\pi)\dot \pi  = 0.$$
 The adjoint of this equation is then
 $$ 0=(\dot \pi)^* (I-\pi)^* = \dot\pi(I-\pi).$$
  Differentiating $\pi^2=\pi$ and applying $\dot\pi=\pi\dot \pi$,
we also have $$\dot \pi \pi =0.$$ Adding these two equations shows
that $$\dot\pi = (I-\pi)\dot \pi + \pi \dot\pi=0$$ in the sense of
distributions. So $\pi$ is constant along almost every coordinate
line segment.  Then it is easy to see that $\pi$ is constant almost
everywhere, and thus is equal to a constant matrix in the sense of
distributions.

We should remark that this simple proof works because $d/dt$ is a
real operator. More properly, on an affine manifold, $\bar
\partial$ is a real operator: We may ignore our convention (\ref{conj-form}),
and instead map $\bar\partial$ to the real operator
$\sfrac12\nabla$ instead via the a natural map from $\mathcal
A^{0,1}(\mbox {End}\, E) \to \Lambda^1(\mbox{End}\,E)$ induced
by $dz^i\mapsto dx^i$. So $\pi^*=\pi$ implies $\dot\pi^* = \dot
\pi$. This fails in the case of complex manifolds, and the
proof to show that the image of $\pi$ is a coherent analytic
subsheaf is quite a bit more involved (Uhlenbeck-Yau
\cite{uhlenbeck-yau86,uhlenbeck-yau89}), although see the
simplification by Popovici \cite{popovici05}.
\end{proof}

\begin{prop}
The flat subbundle $F=\pi(E)\subset E$ is a proper subbundle.  In
other words,
$$ 0 < {\rm rank}\,F < {\rm rank}\,E.$$
\end{prop}
\begin{proof}
First of all, note that ${\rm rank}\,F$ is a constant over $M$,
since it is equal to the rank of $\pi$ as an endomorphism, and $\pi$
is locally constant.

Now $f^0_\infty= \lim_{\sigma\to 0}  f^\sigma_\infty$ is not
identically zero since $f_\infty\neq0$ (Proposition
\ref{construct-f-0-infty}), and the eigenvalues of $f_\infty^\sigma$
are nonnegative and nondecreasing as $\sigma\to0$. So $\pi = I_E -
f^0_\infty$ is not identically $I_E$. Since $\pi$ is a projection,
${\rm rank}\,\pi < {\rm rank}\, E$.

On the other hand, Lemma \ref{basic-rho-lemma} (there is everywhere
on $M$ an eigenvalue of $f_i$ which is bounded by $\rho_i\to0$)
shows that $f_\infty$ has a nontrivial kernel at almost every point.
Therefore, $f^0_\infty$ does as well, and $\pi = I_E-f^0_\infty$
cannot be identically 0. So ${\rm rank}\,\pi >0$.
\end{proof}

\begin{prop}
The flat subbundle $F=\pi(E)$ is a destabilizing subbundle of
$E$. In other words, $$ \frac{{\rm deg}_g E} {{\rm rank}\,E} =
\mu_gE \le \mu_g F = \frac{{\rm deg}_g F} {{\rm rank}\,F}.$$
\end{prop}

\begin{proof}
Recall $$\mu_gE = \frac1r\int_M \frac{c_1(E,h)\wedge \omega_g
^{n-1}} \nu = \frac1{nr} \int_M \mbox{tr}\, K_0 \, \frac{\omega_g
^n} \nu,$$ and for $s = \mbox{rank}\,F$ and $K_F$ the extended mean
curvature of the extended Hermitian connection on $F$ with respect
to the Hermitian metric $h_0|_F$ the restriction of $h_0$ to $F$.
$$\mu_g F = \frac1s\int_M \frac{c_1(F,h_0|_F)\wedge \omega_g
^{n-1}} \nu = \frac1{ns} \int_M \mbox{tr}\, K_F \, \frac{\omega_g
^n} \nu.$$

The Chern-Weil formula (see e.g.\ Kobayashi \cite{kobayashi87})
shows that $\mbox{tr}\, K_F = \mbox{tr}(K_0\pi) - |\pi^\perp
\partial_0 \pi|^2$ for $\pi^\perp
\partial_0 \pi$ the second fundamental form of the subbundle
$F\subset E$.  Now $$\pi^\perp \partial_0 \pi = (I_E - \pi)
\partial_0 \pi = \partial_0\pi - \pi \partial_0\pi =
\partial_0\pi.$$  If we define $K^0 = K_0 - \gamma I_E$, then
$\mbox{tr}\, K^0 = 0$ and $$\mu_g F = \frac1{ns} \int_M
[\mbox{tr}\, (K^0\pi) - |\partial_0\pi|^2] \, \frac{\omega_g
^n} \nu + \frac \gamma n \int_M \frac{\omega_g^n}\nu,$$ while
(\ref{gamma-def}) shows  $\mu_g E = \frac\gamma n \int_M
\frac{\omega_g^n} \nu $.  Therefore, in order to show $\mu_g F
\ge \mu_g E$, we need to show
\begin{equation} \label{int-destabilize}
 \int_M \mbox{tr} (K^0\pi) \frac{\omega_g^n} \nu \ge \int_M |
 \partial_0 \pi |^2 \frac{\omega_g^n} \nu.
\end{equation}

Since $\D\pi = \lim_{\sigma\to0} \lim_{i\to\infty} (I_E -
f^\sigma_i)$ strongly in $L^2$ and $\mbox{tr} \,K^0=0$,
$$\int_M \mbox{tr}(K^0\pi) \frac{\omega_g^n} \nu = -\lim_{\sigma
\to 0} \lim_{i\to\infty} \int_M \mbox{tr}(K^0 f^\sigma_i)
\frac{\omega_g^n} \nu.$$ Compute, using equation
(\ref{cont-method-def}),
\begin{eqnarray*}
-\int_M \mbox{tr}(K^0 f^\sigma_i) \frac{\omega_g^n} \nu &=& \int_M
\epsilon_i \,\mbox{tr}(\log f_{\epsilon_i} \cdot f^\sigma_i)
\frac{\omega_g^n} \nu \\
&& {}+ \int_M \mbox{tr} \{[\mbox{tr}_g \bar\partial (f_i^{-1}
\partial_0 f_i) ] f_i^\sigma \} \frac {\omega_g^n} \nu \\
&\ge & \int_M \mbox{tr} \{ [\mbox{tr}_g \bar\partial (f_i^{-1}
\partial_0 f_i) ] f_i^\sigma \} \frac {\omega_g^n} \nu \\
&=& n\int_M \frac{\mbox{tr} \{[\bar\partial (f_i^{-1}
\partial_0 f_i) ] f_i^\sigma \} \wedge \omega_g^{n-1}} \nu \\
&=& n\int_M \frac{\mbox{tr} [ (f_i^{-1}\partial_0 f_i) \wedge \bar
\partial (f_i^\sigma)] \wedge \omega_g^{n-1} } \nu \\
&& {}+n \int_M \frac{\mbox{tr} [ (f_i^{-1}\partial_0 f_i)
 f_i^\sigma] \wedge \bar\partial \omega_g^{n-1} } \nu, \\
\end{eqnarray*}
where the inequality follows from a local calculation as in
\cite[p.\ 89]{lubke-teleman95} and the last equality follows from
Proposition \ref{int-by-parts} and integration by parts.  Now a
local computation shows that the last integral above satisfies
$$\int_M \frac{\mbox{tr} [ (f_i^{-1}\partial_0 f_i)
 f_i^\sigma] \wedge \bar\partial \omega_g^{n-1} } \nu =
\frac1\sigma \int_M \frac{ \partial [\mbox{tr}(f^\sigma_i) ]\wedge
\bar \partial \omega^{n-1}_g}\nu = 0$$ by integration by parts since
$g$ is affine Gauduchon. On the other hand, other term
 \begin{eqnarray*}
n\int_M \frac{\mbox{tr} [ (f_i^{-1}\partial_0 f_i) \wedge \bar
\partial (f_i^\sigma)] \wedge \omega_g^{n-1} } \nu &=& \int_M
\mbox{tr}\,\mbox{tr}_g [ (f_i^{-1}\partial_0 f_i) \wedge \bar
\partial (f_i^\sigma)] \frac{\omega_g^n} \nu\\
&=&\int_M \mbox{tr}_g\, h_0 (f_i^{-1} \partial_0 f_i,
\partial_0(f_i^\sigma)) \frac{\omega_g^n}\nu \\
&\ge& \int_M |f_i^{-\frac\sigma2} \partial_0( f_i^\sigma) |^2
\frac{\omega_g^n}\nu \\
&\ge& \| \partial_0 (f_i^\sigma)\|^2_{L^2} \\
&=& \| \partial_0 (I_E  - f_i^\sigma) \|^2_{L^2}.
 \end{eqnarray*}
Here, the second line follows from the first since $h_0(A,B) =
\mbox{tr}(AB^*)$ for $B^*$ the $h_0$-adjoint of $B$, the third
line follows by a local computation \cite[Lemma
3.4.4.i]{lubke-teleman95}, and the fourth line follows since
$f_i\le I_E$.

Therefore, $$-\int_M \mbox{tr}(K^0 f_i^\sigma) \frac{\omega_g^n} \nu
\ge \| \partial_0(I_E - f_i^\sigma)\|^2_{L^2},$$ and since
$\partial_0\pi$ is the weak $L^2$ limit of $\partial_0(I_E -
f_i^\sigma)$, $$\lim_{\sigma\to0} \lim_{i\to\infty} \|
\partial_0 (I_E - f_i^\sigma)\| ^2_{L^2} \ge
\| \partial_0 \pi \|^2_{L^2}.$$ This proves the proposition.
\end{proof}

This proposition completes the proof of Theorem \ref{main-thm}.

\section{Simple bundles} \label{simple-section}

Some of this section is a simplified version of Kobayashi
\cite[Section V.7]{kobayashi87}.
\begin{prop} \label{co-simple}
Every $\co$-stable flat vector bundle $E$ over a compact
special affine manifold $M$ is $\co$-simple.
\end{prop}
\begin{proof}
Consider a locally constant section $f$ of $E^*\otimes E$, and let
$a\in\co$ be an eigenvalue of $E$ at a point $x\in M$. Then $f-a
I_E$ is a locally constant endomorphism of $E$ which has a 0
eigenvalue at $x$. Consider $H = (f-a I_E)(E)$. Thus $\mbox{rank}
H<\mbox{rank} E$. We use the $\co$-stability to show $H=0$. If
$\mbox{rank} H>0$, then the stability of $E$ implies that $$\mu(H) <
\mu(E).$$ But we can also identify $H$ with the quotient bundle $E /
\ker (f-aI_E)$, which implies
$$\mu(E)<\mu(H),$$ which provides a contradiction.  Thus $H=0$ and
$f=a I_E$ for the constant $a\in\co$.

Thus the proposition follows from the following
\end{proof}
\begin{prop}
If $E$ is a  $\co$-stable flat vector bundle over a compact
special affine manifold $M$, then any flat quotient vector
bundle $H$ over $E$ satisfies $\mu(E)>\mu(H)$.
\end{prop}
\begin{proof}
If $$0\to F\to E\to H\to0$$ is an exact sequence of flat vector
bundles on $M$, then
 \begin{equation} \label{degree-add}
 \deg F + \deg H = \deg E
 \end{equation}
The proof of (\ref{degree-add}) is to compute the affine first Chern
form.

In terms of a locally constant frame $s_1,\dots,s_r$ of $E$, and for
$h_{\alpha\bar\beta} = h(s_\alpha, s_\beta)$ as above, the first
Chern form is
 \begin{equation}\label{first-chern-formula}
 c_1(E,h) = -\partial \bar\partial \log \det
h_{\alpha \bar \beta}. \end{equation} We will show that there are
natural frames and metrics so that $c_1(E) = c_1(F) + c_1(H)$.

On each sufficiently small open set $U\subset M$, there is a locally
constant frame $\{s_1,\dots,s_r\}$ so that $\{s_1,\dots,s_{r'}\}$ is
a locally constant frame of the subbundle $F$ (for $r'\le r$ the
rank of $F$).  Then the equivalence classes $\{[s_{r'+1}], \dots
[s_r]\}$ form a locally constant frame of the quotient bundle $H$
(here, at $x\in U$, $[s(x)] = s(x) + F_x\in E_x/F_x=H_x$).

We assume $E$ admits a Hermitian metric $h$. Then $h|_F$ is a
Hermitian metric on $F$.  Now there is an orthonormal frame
$\{t_1,\dots, t_r\}$ of $E$ so that $t_1,\dots,t_{r'}$ are sections
of $F$. Then the change-of-frame matrix $A=(A_\alpha^\beta)$
satisfying $t_\alpha = A_\alpha^\beta s_\beta$ is block-triangular
of the form
\begin{equation} \label{block-diag} A =\left(\begin{array}{cc}P & * \\ 0 &Q
\end{array}\right),
\end{equation}
where $P$ is the change-of-frame matrix on $F$ taking $\{s_1,\dots,
s_{r'}\}$ to $\{t_1,\dots,t_{r'}\}$.  The metric $h$ allows us to
identify the quotient bundle $H$ with the orthogonal complement
$F^\perp$ of $F$ in $E$ by orthogonal projection.  Under this
identification, the matrix $Q$ is the change-of-frame matrix on
$F^\perp$ taking $\{[s_{r'+1}], \dots, [s_r]\}$ to
$\{t_{r'+1},\dots, t_r\}$.  Note (\ref{block-diag}) shows $\det A =
(\det P) (\det Q)$.

Now note that the metric $h=(h_{\alpha\bar\beta})$ can be recovered
from a change of frame matrix $A$ by $h=(A\bar
A^\perp)^{-1}$---i.e., $A_\alpha^\gamma h_{\gamma\epsilon} \bar
A_\beta^\epsilon = \delta_{\alpha\beta}$ for the Kronecker
$\delta_{\alpha\beta}$. Then the formulas
(\ref{first-chern-formula}) and (\ref{block-diag}) show that $c_1(E)
= c_1(F) + c_1(H)$.

So the degree addition formula (\ref{degree-add}) follows from the
definition (\ref{degree-def}). Now $$\mu_g(F) < \mu_g(E) \qquad
\Longleftrightarrow \qquad \mu_g(H)> \mu_g(E),$$ which proves the
proposition.
\end{proof}

Finally, we consider the case of real flat vector bundles. Now
let $E$ be a real flat vector bundle over a compact special
affine manifold $M$ equipped with an affine Gauduchon metric
$g$. Such a vector bundle $E$ is said to be $\re$-stable if
every real flat subbundle $F$ of $E$ satisfies
$$ 0 < \mbox{rank}\, F < \mbox{rank}\, E \qquad \Longrightarrow
\qquad \mu_g(F) < \mu_g(E).$$ It is obvious that the
$\co$-stability of $E\otimes_\re\co$ implies the
$\re$-stability of $E$, but the converse may not be true.

\begin{prop}
Let $E$ be an $\re$-stable flat real vector bundle over $M$ a
compact special affine manifold.  As a complex flat vector
bundle, $E\otimes_\re\co$ satisfies one of the following:
\begin{itemize}
 \item $E\otimes_\re\co$ is $\co$-simple.
 \item $E\otimes_\re\co = V \oplus \bar V$, where $V$ is a
     $\co$-stable flat complex vector subbundle of
     $E\otimes_\re\co$ and $\bar V$ is its complex
 conjugate as a subbundle of $E\otimes_\re \co$.
\end{itemize}
\end{prop}
\begin{proof}
Case 1: Every real locally constant section of $\mbox{End}\,E$
has only real eigenvalues at every point $x\in M$. In this
case, let $f$ be a real locally constant section of
$\mbox{End}\,E$, and let $a\in\re$ be an eigenvalue of $f$ at a
point $x\in M$. Then $f-a I_E$ is a real section of
$\mbox{End}\,E$ and, following the proof of Proposition
\ref{co-simple} above, $f-aI_E$ must be identically 0, since
$E$ is $\re$-stable. So $f=a I_E$. The same is true for a
complex locally constant section of $\mbox{End}\,E$ by
considering real and imaginary parts. Thus $E\otimes_\re \co$
is $\co$-simple in this case.

Case 2: There is a real locally constant section $f$ of
$\mbox{End}\,E$ with an eigenvalue $a\notin\re$ at a point
$x\in M$. Then $g=(f-aI)\circ(f-\bar aI)$ is a real section of
$\mbox{End}(E)$. Again, as in the proof of Proposition
\ref{co-simple}, $g$ must be identically 0. So we have the
following splitting into eigenbundles $$E\otimes_\re\co = E_a
\oplus E_{\bar a} = E_a \oplus\overline{E_a}.$$

Now we show that $E_a$ and $E_{\bar a}$ must each be
$\co$-stable. Let $F$ be a flat complex subbundle of $E_a$.
Then it is easy to see that $F\oplus \bar F$ is a real
subbundle of $E_a \oplus \overline{E_a} = E\otimes_\re\co$. The
$\co$-stability of $E_a$ follows from the observation that the
slope $\mu (F) = \mu(F\oplus\bar F)$ for any flat subbundle $F$
of $E_a$.

This observation may be proved by noting that
$\mbox{rank}(F\oplus\bar F) = 2\,\mbox{rank}\,F$, and that the
degree $\mbox{deg}(F\oplus\bar F) = 2\,\mbox{deg}\,F$ also. The
degree calculation can be verified by choosing a Hermitian
metric $h$ on $F$ and extending it to $F\oplus\bar F$ by
setting
\begin{equation} \label{extend-to-conj}
h(\xi, \bar \eta) = h(\bar \xi, \eta) = 0, \qquad h(\bar\xi, \bar
\eta) = \overline{h(\xi,\eta)}
 \end{equation}
for $\xi,\eta$ sections of $F$.
\end{proof}

\begin{cor} \label{re-cor}
Any $\re$-stable flat real vector bundle $E$ over a compact
special affine manifold $M$ admits a real Hermitian-Einstein
metric.
\end{cor}
\begin{proof}
If $E$ is $\co$-stable, then we are done. If not, the previous
proposition shows that $E\otimes_\re\co = V\oplus \bar V$ for
$V$ a complex stable flat subbundle.  Then $V$ admits a
Hermitian-Einstein metric.  It extends to a real
Hermitian-Einstein metric on $E\otimes_\re\co$ by using
(\ref{extend-to-conj}) above.
\end{proof}

\bibliographystyle{abbrv}
\bibliography{thesis}

\end{document}